\documentclass[11pt,letterpaper]{amsart}
\usepackage[english]{babel}
\usepackage[utf8]{inputenc}
\usepackage{lmodern}
\usepackage[T1]{fontenc}
\usepackage{graphicx}
\usepackage{amsmath,amsthm,amssymb}
\usepackage{amsfonts}
\usepackage{amssymb}
\usepackage{enumitem}
\usepackage{stmaryrd}
\usepackage{esint}
\usepackage{soul}
\usepackage{floatrow}
\usepackage{subeqnarray}
\usepackage[all]{xy}
\usepackage{varioref}
\usepackage{tikz-cd}
\usepackage{hyperref}
\usepackage{pgfplots}

\newcommand{\cat}{\operatorname{CAT}}
\newcommand{\ch}{\operatorname{ch}}
\newcommand{\chinf}{\operatorname{ch}(X^\infty)}

\newcommand{\bnd}{\operatorname{bnd}}

\newcommand{\circum}{\operatorname{circ}}
\newcommand{\bdinf}{X^\infty}
\newcommand{\Aut}{\operatorname{Aut}}

\DeclareMathOperator{\supp}{\operatorname{supp}}

\DeclareMathOperator{\iso}{Isom}
\newcommand{\R}{\mathbb{R}}

\newcommand{\proj}{\operatorname{proj}}
\newcommand{\germ}{\operatorname{germ}}

\newcommand{\SL}{\operatorname{SL}}

\newcommand{\PGL}{\operatorname{PGL}}
\newcommand{\res}{\operatorname{Res}}

\newcommand{\bd}{\partial_\infty}

\newcommand{\prob}{\operatorname{Prob}}
\newcommand{\bary}{\operatorname{bar}}
\newcommand{\dt}{\partial T}
\newcommand{\bdg}{\partial_{\operatorname{Grom}}}

\newcommand{\adt}{\tilde{A}_2}

\theoremstyle{plain}
\newtheorem{thm}{Theorem}[section]

\newtheorem{cor}[thm]{Corollary}
\newtheorem{lem}[thm]{Lemma}
\newtheorem{prop}[thm]{Proposition}
\newtheorem{thmi}{Theorem} 

\newtheorem{propi}[thmi]{Proposition}
\theoremstyle{definition} 
\newtheorem{Def}[thm]{Definition}
\theoremstyle{remark}
\newtheorem{rem}[thm]{Remark}

\newtheorem*{ackn}{Acknowledgements}

\title{Stationary measures and random walks on $\tilde{A}_2$-buildings}
\author{Corentin Le Bars}
\date{}

\begin{document}
	\begin{abstract}
		We consider a non-elementary group action $G \curvearrowright X$ of a locally compact second countable group $G$ on a possibly exotic non-discrete affine building $X$ of type $\tilde{A}_2$. We prove that if $\mu$ is an admissible symmetric probability measure on $G$, there is a unique $\mu$-stationary measure supported on the chambers at infinity of the spherical building at infinity. We use this result to study random walks induced by the $G$-action, and we prove that if $\mu$ has finite second moment, $(Z_n o)$ converges almost surely to a regular point of the boundary and the Lyapunov spectrum of the random walk is simple. Applied to Bruhat-Tits buildings, these results extend some classical theorems due to H.~Furstenberg. 
	\end{abstract}
	\maketitle
	
\section{Introduction}
\subsection{Random walks and buildings}
	Tits showed \cite{tits86} that in dimension at least 3, any irreducible affine building is the Bruhat-Tits building of an algebraic group over a field with valuation. Dimension-1 affine buildings correspond to trees or, more generally, $\R$-trees. In dimension 2, there exist irreducible affine buildings that are not Bruhat-Tits: they are said \emph{exotic}. The automorphism groups of these buildings, when non-trivial, are typically discrete \cite{radu19}, and can be cocompact. The simplest and most studied family of exotic buildings are those of type $\tilde{A}_2$: Ronan showed in a constructive way that there is a uncountable family of exotic buildings of type $\tilde{A}_2$ \cite{ronan86}. Discrete buildings of type $\adt$ correspond to buildings in which the apartments are tiled by equilateral triangles: an algebraic example is the Bruhat-Tits building associated to $\SL(3,K)$, for $K$ a discretely valued field.

	Many structural results are known for lattices in higher rank simple algebraic group: Margulis super-rigidity and the normal subgroup theorem are emblematic examples. It is natural to try and extend these results to groups that act by automorphisms on exotic buildings of dimension 2. This field of research has known much progress recently. For instance it was proven in \cite{delaSalle_lecureux_witzel23} that $\adt$-lattices have strong property (T), and in \cite{bader_furman_lecureux23} that such groups also satisfy a normal subgroup theorem. This program implies a translation of methods and results on algebraic groups to a purely geometric and analytic language, and such a translation is interesting and fruitful in itself. In this paper, we consider an action on a $\adt$-affine building $X$, and we investigate questions related to measurable boundary maps and random walks on $X$, in the spirit of Furstenberg's work for linear groups \cite{furstenberg63,furstenberg73}. 
	
	A key feature of Furstenberg's approach to these dynamical systems is the study of stationary measures, i.e. measures that are ``invariant'' under the dynamical system.
	Finding stationary or invariant measures on $G$-spaces and studying their properties is one of the major tasks in the theory of random dynamical systems. 
	
	The goal of this paper is to undertake this when $X$ is a non-discrete affine building of type $\tilde{A}_2$. Having in mind the deep connections between buildings and algebraic groups, studying stationary measures and random processes on $\tilde{A}_2$-buildings can be seen as an extension of the theory of random products in linear groups.

	\subsection{Statements of the results for abstract buildings}

	In the seminal article \cite{furstenberg63}, H.~Furstenberg defines the Poisson-Furstenberg boundary $B $ associated to a probability measure $\mu$ on a group $G$: it is a Borel probability space that represents the future asymptotic directions of the $\mu$-random walk on $G$. Among other features, $B$ can be used to describe all the bounded $\mu$-harmonic functions on $G $ in a minimal way. It is naturally endowed with a non-singular $G$-action which is amenable in Zimmer's sense \cite{zimmer78} and satisfies a number of interesting ergodic properties, see for instance \cite{kaimanovich03}. The notion of a \emph{$G$-boundary} was introduced by U.~Bader and A.~Furman \cite{bader_furman14} to characterize $G$-spaces sharing these properties. An emblematic example of $G$-boundary is indeed the Poisson-Furstenberg boundary associated to an admissible symmetric probability measure on $G$, but they can also arise in other situations. 
	
	Our first result classifies the $G$-maps from any $G$-boundary to the set $\chinf$ of chambers at infinity of a possibly exotic $\mathbb{R}$-building $X$ of type $\tilde{A}_2$. We say that an isometric action $G\curvearrowright X$ on an affine building is \emph{non-elementary} if there is no bounded orbit in $X$, nor finite orbit on the visual boundary $\bd X$.  
	\begin{thmi}[{Boundary maps}]\label{thm intro bd map}
		Let $X$ be a separable complete building of type $\tilde{A}_2$ and let $G$ be a locally compact second countable group. Let $G\curvearrowright X$ be a type-preserving and non-elementary action by isometries. Let $(B, \nu)$ be a $G$-boundary. Then there exists a unique measurable map \[B \rightarrow \ch (\bdinf)\] which is $G$-equivariant. 
	\end{thmi}
	The separability condition is so that the boundary $\bd X$ is metrizable, which is useful for measure-theoretic purposes. Nevertheless, we can get rid of this assumption using a reduction process explained in Section \ref{section reduction intro} below. 
	As this result holds for any $G$-boundary, it can be applied to random walks as well as in other contexts, for instance in order to study more general random dynamical systems like ergodic cocycles, see \cite[Section 5]{bader_furman14}. 
	
	\begin{rem}
		We make the implicit assumption that the action on the building is by automorphisms in the sense of, for instance, \cite[\S 2.1.13]{rousseau23}, that is, isometries that preserve the structure of walls. The group of automorphisms contains a finite-index subgroup formed by automorphisms that are type-preserving, see \cite[\S2.4.6.1]{rousseau23}. Hereafter, we therefore always implicitly make the mild assumption that the $G$-action is type-preserving. Moreover, by \cite[Lemma 3.5]{schillewaert_struyve_thomas22}, we may assume that $X$ is an $\mathbb{R}$-building in which each vertex is special.
	\end{rem}

	Let $\mu$ be a probability measure on $G$. We now apply Theorem~\ref{thm intro bd map} to the study of $\mu$-generated random walk $(Z_n o)$, where $o \in X$ is a basepoint. In the following, we denote as usual by $(\Omega, \mathbb{P}) := (G^{\mathbb{N}}, \mu^{\otimes \mathbb{N}})$ the space of increments associated to $\mu$. By construction, there is a natural measurable factor 
	\begin{align}
		\bnd : (\Omega, \mathbb{P}) \rightarrow (B, \nu_B), \label{eq bnd}
	\end{align}
	from the space of increments to the Poisson-Furstenberg boundary of $(G,\mu)$, which corresponds to the fact that $(B, \nu_B)$ describes the asymptotic behavior of $(Z_n)$. A fundamental result in boundary theory is that to every $\mu$-stationary measure $\nu$ on a standard Borel $G$-space $Y$, one can associate a $G$-equivariant measurable map 
	$$f : b \in B \rightarrow f(b) \in \prob(Y) $$
	from the Poisson-Furstenberg boundary $(B, \nu_B)$ of $(G, \mu)$ to the space of probability measures on $Y$. 
	Conversely, any such map gives rise to a stationary measure defined by 
	$$\nu = \int f(b) d \nu_B(b).$$
	
	Thanks to this duality, we obtain our main result. 
	
	\begin{thmi}[{Stationary measures}]\label{thm intro stat measure}
		Keep the same assumptions on $G, X$ and on the action $G \curvearrowright X$ as in Theorem \ref{thm intro bd map}. Let $\mu $ be an admissible symmetric measure on $G$, and let $(B, \nu_B)$ be the Poisson-Furstenberg boundary of $(G, \mu)$. Then there exists a unique $\mu$-stationary measure $\nu$ on $\ch(\bdinf)$. Moreover, we have the decomposition 
		\begin{eqnarray}
			\nu = \int_{b \in B} \delta_{\psi(b)} d\nu_B(b) \nonumber, 
		\end{eqnarray}
		where $\psi$ is the unique boundary map given by Theorem \ref{thm intro bd map}.
	\end{thmi}

	Two chambers at infinity are said to be \emph{opposite} if they belong to the boundary of a unique apartment. We prove that for the unique stationary measure $\nu $ on $\ch(\bdinf)$ as above, almost every pair of chambers at infinity are opposite. 
	\begin{propi}\label{prop intro opp}
		Let $\nu$ be the unique $\mu$-stationary measure on $\ch(\bdinf)$ given by Theorem \ref{thm intro stat measure}. Then $\nu \otimes \nu$-almost every pair of chambers in $\ch(\bdinf)$ are opposite. 
	\end{propi}
	
	Actually, we prove in Proposition \ref{prop bd map antipodal} that for any $G$-boundary $(B, \nu_B)$, the unique measurable equivariant map $B \to \chinf$ given by Theorem~\ref{thm intro bd map} is almost surely antipodal. Such a statement can be useful in a number of applications. In particular, it is needed for the study of the limit laws of the random walk on $X$ induced by $\mu$. Let $o \in X$ be a basepoint. We say that the measure $\mu$ has \emph{finite second moment} if 
	$$ \int_G d(o, go)^2 d\mu(g) < \infty. $$
	Let $\mathfrak{a}^{+}$ be a fundamental Weyl chamber, and denote by $\mathfrak{a}^{++}$ its interior. For $\lambda \in \mathfrak{a}^{+}$, we say that a sequence $(x_n) $ in $X$ is \emph{$\lambda $-regular} if it sublinearly tracks a $\lambda$-ray, see Section \ref{section reg} for precise definitions. 
	\begin{thmi}[{Simplicity of the Lyapunov spectrum}]\label{thm intro rw}
		Keep the same assumptions on $G, X$ and on the action $G \curvearrowright X$ as in Theorem \ref{thm intro bd map}. Let $\mu$ be a symmetric and admissible probability measure with finite second moment. Then there exists a regular vector $\lambda \in \mathfrak{a}^{++}$ such that the $\mu$-generated random walk $(Z_n o)$ is $\mathbb{P}$-almost surely $\lambda$-regular. In particular, $(Z_n (\omega) o)$ converges $\mathbb{P}$-almost surely to a regular point of the visual boundary contained in the interior of the chamber at infinity $\psi(\bnd(\omega))$, where $\psi$ is the boundary map given by Theorem \ref{thm intro bd map}, and $\mathbb{P}$-almost surely:
		$$ \lim_n \frac{1}{n} \theta(o, Z_n o) = \lambda.$$
	\end{thmi}
	Here we denoted by $\bnd $ the natural projection~\eqref{eq bnd} and by $\theta$ the type function on $X$, see Section \ref{section intro gen imm}. We point out that the moment assumption in the previous theorem is not optimal.  We conjecture that Theorem \ref{thm intro rw} holds if $\mu$ is only assumed to have finite first moment 
	$$ \int_G d(go, o)d\mu(g) < \infty,$$ 
	and that under a finite	second moment assumption, one can prove a central limit theorem using the strategy of Benoist and Quint drawn in \cite{benoist_quint16CLTlineargroups}. Let us also mention that this result relies on the specific geometry of $\adt$-buildings and the fact that we used a symmetric measure: Theorem~\ref{thm intro rw} no longer holds for $\tilde{C}_2$ or $\tilde{G}_2$-buildings even in the algebraic case, as there exist Zariski-dense discrete subgroups of algebraic groups over non-Archimedean local fields associated to these buildings whose limit cones lie in proper facets of the model Weyl chamber \cite{quint02}.
	
	We emphasize that the buildings we consider here are possibly \emph{non-discrete}, which can be interesting for various reasons. As B.~Kleiner and B.~Leeb showed in \cite[Theorem 5.2.1]{kleiner_leeb97}, non-discrete affine buildings arise naturally as asymptotic cones of symmetric spaces of non-compact type and of (even locally finite) affine buildings. The use of non-discrete affine buildings is crucial to their proof of quasi-isometric rigidity of higher-rank symmetric spaces of non-compact type. Another source of interesting examples is that of reductive algebraic groups $\mathbf{G}(K)$ over fields with a dense valuation $v : K \to \mathbb{R}$. This is actually the framework of F.~Bruhat and J.~Tits, who defined affine buildings associated to valuations on root group data which need not be discrete \cite{bruhat_tits72}. Moreover, actions on non-discrete buildings appear on the boundary of the real spectrum compactification of character varieties \cite{burger_iozzi_parreau_pozzetti21}, similar to how actions on $\mathbb{R}$-trees appear in the compactification of the Teichm\"uller space \cite{morgan_shalen84}. Finally, note that while only free groups can act freely (and simplicially) on simplicial trees, some surface groups can act freely on $\mathbb{R}$-trees \cite{morgan_shalen91}. This suggest that actions on non-discrete buildings can be ``richer'' than actions on simplicial ones. 
	
	\subsection{Results for reductive algebraic groups}\label{section intro alg}
	
	The aforementioned results can be translated algebraically in the non-exotic setting. For a general presentation to Bruhat-Tits buildings, we refer to \cite{bruhat_tits72,remy_thuillier_werner15,kaletha_prasad23}. Let $\mathbf{G}$ be a connected reductive isotropic algebraic group defined over a Henselian valued field $K$, and assume that $\mathbf{G}$ admits a valued root group datum as defined in \cite{bruhat_tits72}. Note that the conditions under which one can associate such a valued root group datum to a reductive group over an arbitrary valued field $K$ are subtle if one does not assume that $\mathbf{G}$ is quasi-split over $K$ (or over the maximal \'etale extension of $K$). Some conditions are given in \cite{bruhat_tits84,rousseau23}. Here, we assume that up to performing quasi-split and \'etale descent, there exists a valued root group datum and thus an associated Bruhat-Tits building $X= \Delta^{\mathrm{BT}}(\mathbf{G})$. Note however that this can always be done if the valuation on the Henselian field is discrete with perfect residue field. If the valuation is dense, assume that the field is spherically complete, so that by \cite{martin_schillewaert_steinke_struyve13}, $X$ is metrically complete (this is automatic if the valuation is discrete).  As always, we denote by $G:= \mathbf{G}(K)$ the group of $K$-points of $\mathbf{G}$. The spherical building at infinity $\bdinf= \Delta(\mathbf{G})$ is the Tits building of $(\mathbf{G}, K)$, which can actually be defined over any field. The set of simplices of $\Delta(\mathbf{G})$ (of a given type) is of the form $\mathbf{G}(K)/\mathbf{Q}(K)$, where $\mathbf{Q}(K)$ is (the group of $K$-points of) a parabolic subgroup: in particular we have the identification
	$$\ch(\Delta(\mathbf{G})) \simeq \mathbf{G}(K)/\mathbf{P}(K),$$
	where $\mathbf{P}$ is a minimal $K$-parabolic subgroup of $\mathbf{G}$.

	Under these assumptions, our results translate as the following:

	\begin{thmi}\label{thm intro alg}
		Let $\mathbf{G}$ be an absolutely almost simple algebraic group over a separable and spherically complete field with valuation as above. We assume that the Weyl group of $\mathbf{G}$ is of type $A_2$ (for instance, $G= \SL_3(K)$ or $\PGL_3(K)$). Let $\Gamma$ be a locally compact second countable group and let $\Gamma \to \mathbf{G}$ be a group representation such that the induced $\Gamma$-action on $\Delta^{BT}(\mathbf{G})$ is non-elementary (for instance, $\pi(\Gamma)$ is Zariski-dense and unbounded), and let $(B, \nu_B)$ be a $\Gamma$-boundary. Then there exists a unique measurable $G$-map
		$$B \rightarrow \mathbf{G}(K)/\mathbf{P}(K),$$ for $\mathbf{P}$ a minimal parabolic subgroup of $\mathbf{G}$. 
	\end{thmi}
	
	In particular, if $\Gamma$ is a lattice in a semisimple $\mathbb{K}$-algebraic group  $\mathbf{H}$ for which $\mathbf{P}_\mathbf{H}$ is a minimal parabolic subgroup, and if $\pi : \Gamma \to \mathbf{G}$ is a group homomorphism with Zariski-dense unbounded image, then $ \mathbf{H}(\mathbb{K})/\mathbf{P}_{\mathbf{H}}(\mathbb{K})$ is a $\Gamma$-boundary (a proof of which follows from \cite[Corollary 6.7]{bader_gelander17}) and therefore there exists a unique  $\Gamma$-equivariant measurable map 
	$$ \mathbf{H}(\mathbb{K})/\mathbf{P}_{\mathbf{H}}(\mathbb{K}) \rightarrow \mathbf{G}(K)/\mathbf{P}(K).$$
	This result for local fields is essentially due to Zimmer, and it is the first step in his approach of the proof of Margulis super-rigidity \cite[Chapter 5]{zimmer84}. 
	
	We can also translate the statement about the convergence of the random walk in Theorem \ref{thm intro rw}. We denote by $\mathfrak{a}^+$ a Weyl chamber associated to a choice of positive roots, and by 
	$$\kappa : \mathbf{G} \to \mathfrak{a}^+$$ the Cartan projection given by the Cartan decomposition associated to a choice of a maximal compact subgroup. 
	\begin{thmi}
		Consider the same assumptions as in Theorem \ref{thm intro alg}. If $\mu$ is an admissible symmetric measure on $\Gamma$, then there exists a unique $\mu$-stationary probability measure on the Furstenberg boundary $\mathbf{G}(K)/\mathbf{P}(K)$. If moreover $\Gamma$ is discrete countable and $\mu$ has finite second moment, the Lyapunov spectrum of the random walk $(Z_n (\omega))$ is simple: there exists $\lambda \in \mathfrak{a}^{++} $ such that almost surely,
		$$\lim_n \frac{1}{n}\kappa(Z_n(\omega)) = \lambda. $$
	\end{thmi}
	It should be noted that in the non-Archimedean case, the classification of stationary measures on the flag variety $\mathbf{G}(K)/\mathbf{P}(K)$ is more subtle than in the real case. For instance, for $K=\mathbb{Q}_p$, $p$ a prime number, there exist measures $\mu$ such that the semigroup $\Gamma_\mu$ generated by $\mu$ is Zariski-dense, and for which there exist more than one $\mu$-stationary measures on the flag variety. It is the case for instance if $\Gamma_\mu$ is a small open compact subgroup of $ \mathbf{G}(K)$, see \cite[Remark 4.18]{benoist_quint14}. Here however, the representation is also unbounded, and the uniqueness points to the fact that the action is actually proximal: using the results here, we prove in \cite{le-bars_lecureux_schillewaert23} that the acting group always admit strongly regular hyperbolic elements. 
	
	\subsection{Reduction when the building is not complete separable}\label{section reduction intro}
	
	We end this introduction by discussing the case where the building $X$ is no longer assumed separable and complete. In \cite{le-bars_lecureux_schillewaert23}, we apply a reduction process so that the aforementioned results still hold provided the group $G$ is discrete countable and acts in a non-elementary way. More precisely, we prove that if the $\adt$-building $X$ is not assumed separable nor complete, there exists $Y \subseteq X$ a sub-building of type $\tilde{A}_2$, which is separable and $G$-invariant. The $G$-action on $Y$ is clearly still non-elementary. We can then, using ultralimits of metric spaces, embed $Y$ in a metrically complete $\adt$-building $\hat Y$. The metric completion $\bar Y$ belongs to this ultralimit: it is separable, metrically complete and endowed with a natural $G$-action by isometries. We can then prove that the visual boundary of $\bar Y$ is actually $\bd Y$, therefore the $G$-action on $\bar Y$ remains non-elementary. Although $\bar Y$ is not necessarily a building, the proofs performed here remain true in this setting. Specifically, consider the assumptions: 
	\begin{enumerate}
		\item $X$ is a separable metrically complete $\cat$(0) space, such that $\bd X$ can be given the structure of a spherical building of type $A_2$, and $X$ is contained in a (possibly non-separable) metrically complete $\adt$-building $\hat X$. 
		\item The group $G$ is discrete countable, and acts by isometries on $X$ without bounded orbit in $X $ (equivalently, without fixed point) nor finite orbit in $\bd X$. 
	\end{enumerate}
	As explained before, this situation applies when a discrete countable group $G$ acts on any $\adt$-building $X'$, with $X$ being constructed by the reduction arguments above. Then the following is true:
	\begin{thmi}
		Under the above assumptions, the results stated in Theorems \ref{thm intro bd map}, \ref{thm intro stat measure}, \ref{thm intro rw} and Proposition \ref{prop intro opp} hold. 
	\end{thmi}
	
	In particular, for the algebraic case, it means that our results (when $\Gamma$ is discrete countable) remain true if the non-Archimedean field is only assumed to be complete, and no longer only for fields that are spherically complete and separable. 
	
	\subsection{Structure of the proof}
	An essential ingredient for the existence of the boundary map $B \rightarrow \ch (\bdinf)$ is a result due to Bader, Duchesne and L\'ecureux in \cite{bader_duchesne_lecureux16}. The authors prove that if $B$ is a $G$-boundary, then for a non-elementary $G$-action on a $\cat$(0) space of finite telescopic dimension, there exists a boundary map $\psi: B \rightarrow \bd X$. We investigate the characteristics of this map when $X$ is a building of type $\tilde{A}_2$. If almost surely, $\psi(b)$ is a regular point of the boundary, i.e. is not represented by a panel, we get a $G$-map $B \to \ch (\bdinf)$, which is what we want. If the target space is the set of vertices at infinity $\bdinf_1$, we use the theory of measurable metric fields and a theorem due to B.~Duchesne \cite{duchesne13} generalizing Adams-Ballmann equivariant theorem \cite{adams_ballmann98} in order to exhibit several mutually exclusive situations involving the field of panel trees over $\bdinf_1$. We then use the strong ergodic properties of $G$-boundaries to rule out some of these possibilities, and conclude with the bijection between the ends of panel trees and the residues of vertices a infinity. 
	
	Given a measure $\mu$ on $G$, the boundary map $B(G, \mu) \to \chinf$ produces a stationary measure. Proving uniqueness amounts to showing that any $G$-equivariant map $B(G, \mu) \to \prob(\chinf)$ must factor through this boundary map. 
	
	Finally, we use the notion of $\lambda$-ray developed by V.~Kaimanovich \cite{kaimanovich89} for symmetric spaces and J.~Parkinson and W.~Woess in the context of buildings \cite{parkinson_woess15}. We apply a combination of results by H.~Izeki \cite{izeki23} and A.~Karlsson and G.~Margulis \cite{karlsson_margulis99} to our context, and we conclude with the particular behavior of the opposition involution in $\tilde{A}_2$-buildings. 
	
	\subsection{Structure of the article}
	In Section \ref{section intro immeuble}, we review the definition of (non-discrete) affine building as well as some related objects, including the panel trees associated to vertices at infinity. Section \ref{section bd theory} presents some notions of boundary theory such as measurable metric fields and relative isometric ergodicity. In Section \ref{section meas struc A2}, we ensure that there exist nice measurable structures for some fibered spaces, so that we can apply the results from Section \ref{section bd theory}. We prove Theorem \ref{thm intro bd map} in Section \ref{section bd map immeuble}. In Section \ref{section uniq stat meas A2}, we prove the main Theorem \ref{thm intro stat measure}, and we derive Proposition \ref{prop intro opp}. In Section \ref{section conv}, we show that the random walk associated to an admissible measure with finite second moment on $G$ converges almost surely a regular point of the boundary. 
	\begin{ackn}
		The author is very grateful to J.~L\'ecureux, who originated many ideas implemented here and for helpful discussions. 
	\end{ackn}
	
	\section{Affine buildings of type $\tilde{A}_2$}\label{section intro immeuble}
	
	In this section, we define non-discrete affine buildings, and we introduce the objects that we will use throughout the paper. Standard references are \cite{parreau00} and \cite{rousseau23}. 
	
	\subsection{Non-discrete affine buildings}\label{section intro gen imm}
	
	Let $(W,V)$ be an affine reflection system, that is $W = W_0 \ltimes T$, where $W_0$ is a finite reflection group and $T$ is a translation group on $V=\mathbb{E}^n$. Let $\mathfrak{a}^{+}$ be some Weyl chamber for $W$, and by $\mathfrak{a}^{++}$ its interior. Let $X$ be a set, and let $\mathcal{A}$ be a collection of injective charts of $V$ into $X$, which we call an \emph{atlas}. Each such injection is called a \emph{chart}, or \emph{marked apartment}, and the image $A$ of $V$ by an injection is called an \emph{apartment}. We say that $(X, \mathcal{A})$ is an affine building modelled after $(W, V)$ if the following axioms are verified. 
	\begin{enumerate}[label=(A\arabic*)]
		\item The atlas $\mathcal{A}$ is invariant by pre-composition with $W$. \label{A1}
		\item Given two charts $f, f' : V \to X$ with $f(V) \cap f'(V) \neq \emptyset$, then $U:=f^{-1} (f'(V))$ is a closed convex subset of $V$, and there exists $w \in W$ such that $f|_U = f' \circ w |_{U}$. \label{A2}
		\item For any pair of points $x, y \in X$, there is an apartment containing both. \label{A3} 
	\end{enumerate}
	Axioms \ref{A1}-\ref{A3} imply the existence of a well-defined distance function $d : X \times X \to \R_+ $, such that the distance between any two points is the $d_V$-distance between their pre-image under any chart containing both. The metric space $(X,d)$ is then a CAT(0) space. Every automorphism of $X$ induces an isometry of $(X,d)$.
	A \emph{Weyl chamber} (or \emph{sector}) in $X$ is the image of an affine Weyl chamber under some chart $f \in \mathcal{A}$. 
	\begin{enumerate}[label=(A\arabic*), resume]
		\item Given two Weyl chambers $S_1, S_2$ in $X$, there exist sub-Weyl chambers $S_1'\subseteq S_1$ and $S_2' \subseteq S_2$ such that $S_1'$ and $S_1'$ are contained in the same apartment. \label{A4} 
		\item For any apartment $A$ and $x \in X$, there exists a retraction $\rho_{A, x}: X \to A$ such that $\rho_{A, x} $ does not increase distance and $\rho^{-1}_{A, x} (x ) = \{x\}$. \label{A5} 
	\end{enumerate}
	
	If the affine reflection group is not discrete, we say that the building $(X,\mathcal{A})$ modelled after $(W, V)$ is \emph{non-discrete}. We will assume that the system of apartments $\mathcal{A}$ is maximal. We say that $X$ is of \emph{type $\tilde{A}_2$} if $W_0$ is a spherical Coxeter group of type $A_2$. 
	
	By the axioms \ref{A1}-\ref{A3}, there exists a marked apartment $f$ sending the fundamental closed Weyl chamber $\mathfrak{a}^+$ to a Weyl chamber in $X$ based at $x $ and containing $y$. The \emph{type} $\theta(x,y)$ of the Euclidean segment $[x,y]$ is the unique vector in $\mathfrak{a}^+$ such that $y = f(\theta(x,y))$. We say that the segment $[x,y]$ is \emph{regular} if the type $\theta(x,y)$ is regular, i.e. $\theta(x,y) \in \mathfrak{a}^{++}$. If $[x,y]$ is regular, then both $\theta(x,y)$ and $\theta(y,x)$ are regular. The group $G$ is type-preserving by assumption: this implies that for any $g\in G$, $\theta(gx,gy)=\theta(x,y)$. 
	
	Denote by $w_0$ the long element of the finite Weyl group $W_0$ associated with the affine reflection group $W$. For $\lambda \in \mathfrak{a}$, the \emph{opposition involution} $\iota : \mathfrak{a}^{+} \to \mathfrak{a}^{+}$ is defined by 
	$$ \iota (\lambda) = w_0(-\lambda).$$

	\subsection{Visual boundary and spherical building at infinity}\label{section spherical boundary}
	
	Let $(X,\mathcal{A}, d)$ be an affine building modelled after $(W, V)$. We denote by $\overline{X} = X \cup \bd X$ its visual bordification as a $\cat$(0) space. When $X$ is separable, the visual boundary $\bd X$ is endowed with a natural metrizable topology and isometries of $X$ extend to homeomorphisms on the boundary. 
	
	The boundary $\bd X$ of a $\cat$(0) space $X$ can also be endowed with the Tits metric, with which it becomes a $\cat(1)$ space and on which $\iso(X)$ acts by isometries. The following proposition will be useful later.

	\begin{prop}[{\cite[Proposition 1.4]{balser_lytchak05}}]\label{circumcenter}
		Let $X$ be a $\cat$(1) space of finite dimension and of radius $r \leq \pi/2$. Then $X$ has a circumcenter which is fixed by every isometry of $X$. 
	\end{prop}

	As $X$ is an affine building, we can give $\bd X$ the structure of a spherical building, which we denote by $\bdinf$, see \cite[\S3.2]{rousseau23}. Apartments of $X$ are in bijection with apartments of $\bdinf$. A top-dimensional simplex of $\bdinf$ is a \emph{chamber at infinity}, and two chambers of $\bdinf$ are called \emph{opposite} if there is a unique apartment containing them. For $F$ a face of $\bdinf$ and $x \in X$ a special vertex, we denote by $Q(x,F)$ the affine facet in $X$ with tip $x$ asymptote to $F$, see for instance \cite[Corollary~1.9]{parreau00}. In particular, for $C \in \chinf$ a chamber at infinity, $Q(x,C)$ denotes the Weyl chamber with tip $x$ asymptote to $C$.

	The set of chambers at infinity $\ch (\bdinf)$ can be endowed with a topology, which makes it a totally disconnected space. A basis of open neighborhoods in $\ch (\bdinf)$ is given by 
	\begin{eqnarray}\label{eq basis chinf}
		U_x(y) := \{ C \in \ch(\bdinf) \, | \, y \in Q(x, C)\} \subseteq\chinf, 
	\end{eqnarray}
	for $x,y \in X$. 
	The following proposition summarizes some properties of this bordification. 
	
	\begin{prop}[{\cite[\S 3.2]{rousseau23}}]
		Let $X$ be any Euclidean building. Then there is a topology on $X \cup \chinf$ for which a basis of open sets of the chambers at infinity is given by the sets \eqref{eq basis chinf}. It agrees with the $\cat$(0) topology on $X$. This topology is first-countable and Hausdorff. For every boundary point $\xi \in\bd X$ that is strictly supported on a chamber, one can associate canonically a chamber $C_\xi \in \chinf$. This map is a homeomorphism with respect to the restriction of the visual topology on $\bd X$ to such points. 
	\end{prop}
	A more detailed account on how to give a topology on the set $X \cup \chinf$, and even on $X \cup X^\infty_\tau$, where $X^\infty_\tau$ represents the set of simplices of type $\tau$ of the spherical building at infinity, is given in \cite{rousseau23}. 
	Note in particular that any type-preserving automorphism of $X$ induces a type-preserving automorphism of the spherical building at infinity $X^\infty$. With the topology on $\chinf$, an automorphism of $X$ can be extended to a homeomorphism on $\chinf$. In this paper, any action of a group on an affine building is by automorphisms.

	\subsection{Retractions}\label{section retraction}
	
	Let $X$ be an affine building, $A$ an apartment, $x \in A$ a vertex and let $S$ be a Weyl chamber in the apartment $A$, based at $x$. Then there exists a unique retraction map $\rho_{A, S}: X \to A$ such that $\rho_{A, S}$ preserves the distance on any apartment containing a sub-Weyl chamber of $S$, see \cite[Proposition 1.20]{parreau00}. Moreover, $\rho_{A, S}$ does not increase distances. We call $\rho_{A,S}$ the \emph{(sector-based) retraction} of $X$ on $A$ based at $S$.

	\subsection{Panel trees}\label{section panel trees}
	Let $X$ be an affine building and let $F^\infty$ be a simplex at infinity. There is a way of attaching an affine building $X(F^\infty)$ to the panel $F^\infty$, as a set of equivalence classes of affine panels representing $F^\infty$. If $F^\infty $ is of codimension 1, $T_{F^\infty} := X(F^\infty)$ is an $\mathbb{R}$-tree called the \emph{panel tree at infinity}, see \cite{tits86}. Let $S$ be a sector in $X$, one of whose faces represents $F^\infty$, that is, the parallel class $[S]_\sim$ of $S$ belongs to the residue $\res(F^\infty)$ of $F^\infty$ in the building at infinity $\bdinf$.  We have a well-defined homomorphism of affine buildings 
	\begin{eqnarray}
		\pi_{F^\infty} : X &\longrightarrow& T_{F^\infty} \nonumber \\
		x &\longmapsto& [Q(x, F^\infty)] \nonumber, 
	\end{eqnarray} 
	where  $[Q(x, F^\infty)]$ is the class of the panel based at $x$, in the direction of $F^\infty$. 
	
	Let $S, S'$ be Weyl chambers in $X$ such that $F^\infty \in \partial_\infty S$ and $F^\infty \in \partial_\infty S'$. In other words, the parallel classes of $S$ and $S'$ are in the residue of $F^\infty$. If we assume furthermore that $S$ and $S'$ are parallel, then there is a sub-sector $S'' \subseteq S \cap S'$, and $\pi_{F^\infty}(S'')$ is a common sub-sector of the sectors $\pi_{F^\infty}(S)$ and $\pi_{F^\infty}(S')$ in the affine building $X(F^\infty)$. This allows to define a surjective map 
	\begin{eqnarray}
		\phi_{F^\infty} &:& \res(F^\infty) \longrightarrow \dt_{F^\infty} \nonumber \\
		&&[S]_\sim \longmapsto [\pi_{F^\infty} (S)]_\sim \nonumber,  
	\end{eqnarray}
	where $\dt_{F^\infty}$ is to be understood as the set of chambers at infinity of the affine $\tilde{A}_1$-building $T_{F^\infty}$. 
	We now prove that $\phi_{F^\infty}$ is a homeomorphism, when we endow $\res{(F^\infty)}$ and $\dt_{F^\infty}$ with the natural topology for chambers at infinity on affine buildings. 
	This result is not surprising to experts, but we did not find any reference in the literature in the non-discrete case. In the discrete case, it is proven in \cite[Lemma 4.2]{remy_trojan21}. We denote by $\Aut(X)_{F^\infty}$ the group of automorphisms of $X$ that stabilize $F^\infty$.
	
	\begin{prop}\label{bijection ends chambers}
		Let $ F^\infty$ be a panel in the spherical building at infinity $\bdinf$. Then the map 
		$$\phi_{F^\infty} : \res(F^\infty) \longrightarrow \dt_{F^\infty}$$ 
		is an $\Aut(X)_v $-equivariant homeomorphism between the residue $\res(F^\infty)$ and the space of ends of the panel tree $T_{F^\infty}$.  
	\end{prop}
	
	\begin{proof}
		The fact that $\phi_{F^\infty}$ is a bijection is well-known, see for example \cite[Proposition 4]{tits86} or \cite[Proposition 11.22]{weiss09}. It is $\Aut(X)_{F^\infty}$-equivariant by construction. 
		
		Let $\xi, \eta$ be two points of $T_{F^\infty} \cap \pi_{F^\infty}(X_s)$, and let $\mathcal{C} \in \dt_{F^\infty}$. In order to avoid confusion, we denote by $Q_{F^\infty}(\xi, \mathcal{C})$ the sector (which is a ray since $T_{F^\infty}$ is a tree) based at $\xi$ in the direction of the end $\mathcal{C}$. We shall write $U^\infty_{\xi}(\eta)$ to denote the set 
		\begin{eqnarray}
			\{ \mathcal{C} \in \dt_{F^\infty} \, :  \, \eta \in Q_{F^\infty}(\xi, \mathcal{C})\}  \nonumber. 
		\end{eqnarray}
		
		Let $C \in \phi_{F^\infty}^{-1} (U^\infty_{\xi}(\eta))$ be a chamber at infinity. 
		The sector $Q_{F^\infty}(\xi, \phi_{F^\infty}(C))$ belongs to an apartment $A_\infty$ of the panel-tree $T_{F^\infty}$. 
		Such an apartment can be written $\pi_{F^\infty}(A)$, for $A$ an apartment in $X$ such that $\partial_\infty A $ contains $F^\infty$. 
		Let $x \in X_s \cap A$ be such that $\pi_{F^\infty}(x) = \xi$. 
		By surjectivity, there exists a sector $S$ in $A$ based at $x$ such that $\pi_{F^\infty}(S) $ is a sector of $T_{F^\infty}$ representing $\phi_{F^\infty}(C)$. Let $y $ be a special vertex in $S$ be such that $\pi_{F^\infty}(y) = \eta$. Then by injectivity,
		\begin{eqnarray}
			Q_{F^\infty}(\xi, \phi_{F^\infty}(C)) = \pi_{F^\infty}(S) \nonumber. 
		\end{eqnarray} 
		Now for any $C' \in \res(F^\infty) \cap U_x(y) $, we have $y \in Q(x, C')$ and thus $\eta \in Q_{F^\infty}(\xi, \phi_{F^\infty}(C))$. In other words, $U_x(y) \cap \res(F^\infty) $ is contained in $\phi_{F^\infty}^{-1} (U^\infty_{\xi}(\eta))$, and $\phi_{F^\infty}$ is continuous. 
		
		Conversely, let $x, y \in X$ be special vertices and let $\xi \in \phi_{F^\infty}(U_x(y)\cap \res(F^\infty))$. Take $C \in U_x(y)\cap \res(F^\infty)$ such that $\phi_{F^\infty}(C) = \mathcal{C}$. Define $\xi = \pi_{F^\infty}(x) $ and $\eta = \pi_{F^\infty}(y) $. By definition, $Q_{F^\infty}(\xi, \mathcal{C})$ is a sector of $T_{F^\infty}$. But by injectivity, $Q_{F^\infty}(\xi, \mathcal{C})= \pi_{F^\infty}(Q(x, C))$ and $y \in Q(x, C)$. Therefore, $\eta \in Q_{F^\infty}(\xi, \mathcal{C})$ and $\mathcal{C} \in U^\infty_{\xi}(\eta)$. 
	\end{proof}
	
	Let us now describe these panel trees when $X$ is an affine building of dimension 2. 
	Recall that a boundary point $v \in  \bd X$ is an equivalence class of rays, two rays $r_1$ and $r_2$ being equivalent, for which we write $r_1 \sim r_2$, if they contain subrays that lie in a common apartment and are parallel in this apartment.  We will say that two geodesic rays $r_1$ and $r_2$ are \emph{strongly asymptotic}, and write $r_1 \simeq r_2$, if their intersection contains a geodesic ray.  For two equivalent geodesic rays $r_1$ and $r_2$ that represent the boundary point $v \in \bd X$, we define their distance to be:
	
	\begin{equation*}\label{eq metric panel tree}
		d_v (r_1, r_2) := \underset{s}{\inf} \underset{t \rightarrow \infty}{\lim} d(r_1(t+s), r_2(t)). 
	\end{equation*}
	It defines a pseudo-distance \cite[Section 5.3]{caprace_lecureux11}, and two strongly asymptotic rays  $r_1 $ and $r_2$ satisfy $d_v(r_1, r_2) =0$. This pseudo-distance does not depend on the representatives within the $\simeq $-strongly asymptotic classes: on these $\simeq$-classes, it becomes a distance. The metric space $(T_v, d_v)$ of strongly asymptotic classes of rays in the class of a vertex at infinity is the panel tree at $v$. The branch points of this tree correspond to thick walls of $X$ \cite[\S 4.26]{kramer_weiss14}.
	Note that given a vertex $v \in \bdinf$, there is as before a well-defined and continuous application defined by $\pi_v : X \longrightarrow T_v$, associating to $x \in X$ the $\simeq$-class of the geodesic ray based at $x$ in the direction of $v$. 
	
	Notice that in this case the building $\bdinf$ is endowed with a coloring that divides its vertices into two types, and the corresponding subsets are denoted $\bdinf_1$ and $\bdinf_2$. 
	
	\section{Boundary theory and metric fields}\label{section bd theory}
	
	\subsection{Boundary theory}
	
	Let $G$ be a locally compact second countable group, and let $\mu$ be a probability measure on $G$.
	Let $(\Omega, \mathbb{P}) := (G^{\mathbb{N}}, \mu^{\otimes \mathbb{N}})$ be the \emph{space of forward increments}, with the product $\sigma$-algebra. The application 
	\begin{equation*}
		(\omega=(\omega_i)_{i\in \mathbb{N}},n) \in  \Omega \times\mathbb{N} \mapsto Z_n(\omega) = \omega_0 \omega_1 \dots \omega_{n-1},
	\end{equation*}
	defines the $\mu$-generated (right-)random walk on $G$. 
	We say that $\mu \in \prob(G) $ is \emph{admissible} if its support $\supp(\mu)$ generates $G$ as a semigroup, and there exists $m \in \mathbb{N}$ such that the $m$-th convolution power $\mu^{\ast m}$ and the Haar measure on $G$ are non-singular. We denote by $\check{\mu}= i_\ast \mu$ the inverse probability measure, where $i(g) = g^{-1}$. We say that $\mu \in \prob(G)$ is \emph{symmetric} if $\mu = \check{\mu}$.

	Let $(Y, \mathcal{Y}, \nu)$ be a standard probability space, here defined as a standard Borel space endowed with a probability measure compatible with the Borel structure.  Let $G \curvearrowright Y$ be a Borel action. The convolution probability measure $\mu \ast \nu$ is defined by: for any $f $ bounded measurable function on $Y$, 
	\begin{equation*}
		\int_Y f(y)d(\mu \ast \nu)(y) = \int_G\int_Y f(g \cdot y ) d\mu(g) d\nu (y).
	\end{equation*}
	\begin{Def}
		A probability measure $\nu \in \prob(Y) $ is \emph{$\mu$-stationary} if $\mu \ast \nu = \nu $. In this case, we say that $(Y, \nu)$ is a \emph{$(G, \mu)$-space}. 
	\end{Def}
	The set of probability measures $\prob(Y)$ on $Y$ can be endowed with the metrizable weak-$\ast$ topology. The following statement, due to H.~Furstenberg, is one of the fundamental results of boundary theory, see \cite[Lemme 3.2]{benoist_quint11}. 
	
	\begin{thm}\label{thm lim meas furst}
		Let $(Y, \nu)$ be a $(G, \mu)$-space. Then there exists an essentially well-defined $G$-equivariant measurable map $\omega \mapsto \nu_\omega \in \prob(Y)$ such that $\mathbb{P}$-almost surely, $Z_n (\omega)\nu \rightarrow \nu_\omega $ in the weak-$\ast$ topology. Moreover, we have the decomposition $$\nu = \int_{\Omega} \nu_\omega d \mathbb{P}(\omega). $$
	\end{thm}
	
	A $(G, \mu)$-space $(Y, \nu)$ is a \emph{$(G,\mu)$-boundary} if for $\mathbb{P}$-almost every $\omega\in \Omega$, the measure $\nu_\omega$ given by Theorem \ref{thm lim meas furst} is a Dirac measure. 
	
	For $(Y, \nu)$ and $(Y', \nu')$ non-singular $(G, \mu)$-spaces, we say that a measurable map $\pi :(Y, \nu) \to (Y', \nu')$ is a \emph{factor map} if $\pi_\ast \nu = \nu'$. It is straightforward to see that any measurable $G$-equivariant factor $(Y', \nu')$ of a $(G, \mu)$-boundary $(Y, \nu)$ is still a $(G, \mu)$-boundary. The \emph{Poisson-Furstenberg boundary} $(B(G, \mu), \nu_B)$ associated to $\mu$ is a $(G, \mu)$-boundary that is maximal for this property: for every $(G, \mu)$-boundary $(Y, \nu)$, there 
	is a $G$-equivariant measurable factor map $p : (B, \nu_B) \rightarrow (Y, \nu)$, uniquely defined up to $\nu_B$-null sets. Note that there is a natural measurable projection $\bnd : \Omega \rightarrow B$. 
	
	The same construction can be done for $\check{\mu}$: we denote by $(\check{B}, \check{\nu})$ the Poisson-Furstenberg boundary associated with $\check{\mu}$. Note that $(\check{B}, \check{\nu})$ is a $(G, \check{\mu})$-boundary. 
	
	\subsection{$G$-boundaries}\label{section bd theory immeuble}
	
	Throughout this section, $G$ is a locally compact second countable group. Let $(S, \eta)$ be a non-singular $G$-space and let $S'$ be a standard Borel space endowed with a Borel $G$-action. A map $(S, \eta) \to S'$ between $G$-spaces is called a \emph{$G$-map} if it is a $G$-equivariant measurable map. 
	
	The following notions, which will be crucial in this paper, are taken from \cite{bader_furman14}.
	\begin{Def}[Relative isometric ergodicity]
		Let $(S, \nu)$ and $(Q, \nu')$ be two standard probability $G$-spaces. A $G$-factor $\pi : S \rightarrow Q $ is said to be \emph{relatively isometrically ergodic} (RIE) if for all fiberwise $G$-action on a $G$-map $p : \mathbb{X} \rightarrow A$ for which there exist $G$-equivariant maps $f : S \rightarrow \mathbb{X}$ and $f_0 : Q\rightarrow  A$ such that $p \circ f = f_0 \circ \pi$, then there exists a relative section, i.e. a $G$-map $\phi : Q  \rightarrow \mathbb{X} $ such that $f = \phi \circ \pi$ and $f_0 = p \circ \phi$ almost everywhere. 
	\end{Def}

	\begin{Def}[$G$-boundary pair]
		Let $(B_{-}, \nu_{-})$ and $(B_{+}, \nu_{+})$ be $G$-spaces. We say that $(B_{-}, B_{+} )$ is a $G$-boundary pair if: 
		\begin{itemize}
			\item the actions $G \curvearrowright B_{+}$ and $G \curvearrowright B_{-}$ are Zimmer amenable, see \cite{zimmer84} and \cite[Section 7.1]{duchesne13}; 
			\item the $G$-maps $B_{-} \times B_{+} \rightarrow B_{-} $ and $B_{-} \times B_{+} \rightarrow B_{+} $ are relatively isometrically ergodic. 
		\end{itemize}
		We say that a $G$-space $(B, \nu_B)$ is a $G$-boundary if $(B,B)$ is a $G$-boundary pair. 
	\end{Def}
	
	When $G$ is locally compact second countable, there always exists a $G$-boundary. More precisely, we have the following. 
	
	\begin{thm}[{\cite[Theorem 2.7]{bader_furman14}}]\label{thm poisson bd pair}
		Let $G$ be a locally compact second countable group, and let $\mu$ be an admissible probability measure on $G$, and let $(B, \nu_B)$ be the Poisson-Furstenberg boundary associated with $(G, \mu)$, resp. $(\check{B}, \check{\nu})$ the Poisson-Furstenberg boundary associated with $(G, \check{\mu})$. Then $(B, \check{B})$ is a $G$-boundary pair. 
	\end{thm}

	A reason why it is so convenient to use $G$-boundaries is that they give rigidity results on the group $G$, one of which is the existence of Furstenberg maps. 
	
	\begin{thm}[{\cite[Theorem 34]{bader_duchesne_lecureux16}}]\label{furst map}
		Let $X$ be a complete $\cat$(0) space of finite telescopic dimension and let $G$ be a locally compact second countable group acting continuously by isometries on $X$ without invariant flats. Let $(B_+, \nu_+) $ and $( B_-, \nu_-)$ form a $G$-boundary pair. Then there exist measurable $G$-maps $\varphi_\pm: B_\pm \rightarrow \bd X$. 
	\end{thm}

	\subsection{Metric fields and equivariant theorem}
	
	In the course of the proof, we will use the vocabulary of measurable metric fields, on which we can find accounts in \cite{duchesne13}. A metric field is here a measurable family of metric spaces conveniently attached to a base, and this notion has nothing to do with that of a field with an absolute value. This notion bears many similarities with the language of fiberwise isometric actions, as explained in \cite[Lemma 4.11]{duchesne_lecureux_pozetti23}.
	
	\subsubsection{Measurable metric fields}
	Let $(A, \eta)$ be a measurable space. A measurable metric field $\mathbf{X}$  over $A$ can be thought of as a way of attaching a metric space $X_a$ to each point $a \in A$  in a measurable way. 
	
	\begin{Def}
		Let $(A, \eta)$ be a standard probability space. Let $\mathbf{X}$ be a collection $\mathbf{X} = \{(X_a, d_a)\}_{a \in A}$ of complete separable metric spaces over $(A, \eta)$, and let $\{x^n\}_{n \in \mathbb{N}}$ be a countable set of elements of $\underset{a \in A}{\Pi} X_a $ such that: 
		\begin{itemize}
			\item $\forall n, m \in \mathbb{N}$, the map $ a \mapsto d_a(x^n_a, x^m_a)$ is measurable; 
			\item for $\eta$-almost every $a \in A$, the set $\{x^n_a\}_{n \in \mathbb{N}}$ is dense in $X_a$. 
		\end{itemize}
		The family $\{x^n\}_{n \in \mathbb{N}}$ is called a \emph{fundamental family} for $\mathbf{X}$. The data $((A, \eta), \mathbf{X}, \{x^n\}_{n \in \mathbb{N}})$ is called a \emph{metric field} over $(A, \eta)$. 
		
		A \emph{section} $s : (A, \eta) \rightarrow \mathbf{X} $ is a map such that  $s(a) \in X_a$ for every $ a \in A$ and such that for every element $x^n $ of the fundamental family, $a \mapsto d_a (x^n_a, s(a))$ is measurable. 
	\end{Def}
	
	Let $(A, \eta)$ be a $G$-space, and let $\mathbf{X} = \{X_a\}_{a \in A}$ be a metric field over $A$. We say that $G$ acts on $\mathbf{X}$ by the cocycle $\alpha$ if: 
	\begin{itemize}
		\item for all $g \in G$ and almost every $a \in A$, $\alpha (g, a) \in \iso(X_a, X_{ga})$; 
		\item for all $(g,g') \in G^2$ and almost every $a \in A$, we have $\alpha (gg', a) = \alpha(g, g'a) \alpha(g', a)$; 
		\item for all elements $x^n $ and $ x^m $ of the fundamental family associated to $\mathbf{X}$, the map $(g, a) \mapsto d_a(x^n_a, \alpha (g, g^{-1}a) x^m_{g^{-1}a})$ is measurable. 
	\end{itemize}
	A section $s : (A, \eta) \rightarrow \{X_a\}_{a \in A} $ is then said to be \emph{$\alpha$-invariant} if for all $g \in G$ and almost every $a \in A$, we have $s(ga) = \alpha (g, a) s(a)$. 
	
	There exists a measurable structure on the collection $\prob(\mathbf{X}) := \{\prob(X_a)\}_{a \in A}$, coming from a Borel structure on the set $C(\mathbf{X}):= \{C(X_a)\}_a$  of Banach spaces of continuous functions over $\mathbf{X}$. From now on, we use this structure, and we refer to \cite[Theorem 2.19]{anderegg_henry14} and \cite[Section 9]{duchesne13} for details. 
	
	\subsubsection{$\cat(0)$-fields}
	We say that $\mathbf{X}$ is a \emph{$\cat(\kappa)$-field} if for almost every $a$, $X_a$ is a $\cat(\kappa)$ space. Similarly, if for almost every $a$, $X_a$ is a Euclidean space, we say that $\mathbf{X}$ is a \emph{Euclidean field}. 
	A \emph{subfield} of $\mathbf{X}$ is a collection $\{Y_a\}_a $ of non-empty closed convex subsets such that for every section $x $ of $\mathbf{X}$, the function $a \mapsto d (x_a, Y_a)$ is measurable. A subfield is called \emph{invariant} if for all $g$, and $\eta$-almost every $a \in A$, 
	\begin{eqnarray}
		\alpha(g , g^{-1} a)Y_{g^{-1}a } = Y_a \nonumber. 
	\end{eqnarray}
	
	Let $\mathbf{X}$ be a $\cat$(0)-field over $(A, \eta)$. There is a way to define the a metric field structure on $\partial \mathbf{X} = \sqcup_{a \in A} \bd X_a $ over $\mathbf{X}$ using compactifications by horofunctions, as explained in \cite[Section 8]{duchesne13}. This structure is called the \emph{boundary field} over $\mathbf{X}$.

	The following result about sections of circumcenters will be useful on several occasions. 
	
	\begin{prop}[{\cite[Lemma 8.7]{duchesne13}}]\label{prop section circ borel}
		Let $\mathbf{X}$ be a metric field of complete $\cat(0)$ spaces over a space $(A, \eta)$, and let $\{B_a\}_{a \in A}$ be a Borel subfield of $\mathbf{X}$. Then the circumradius function $a \in A \mapsto r_a(B_a)$ is measurable. If it is essentially bounded, the section of circumcenters 
		\begin{eqnarray}
			\circum : \mathbf{B} &\rightarrow& \mathbf{X} \nonumber \\
			B_a & \mapsto & \circum (B_a), \nonumber
		\end{eqnarray}
		where $\circum (B_a)$ is the circumcenter of $B_a$, is measurable. 
	\end{prop}
	
	\section{Measurable structures and barycenters}\label{section meas struc A2}
	
	Let now $X$ be a complete separable building of type $\tilde{A}_2$ and let $G$ be a locally compact second countable group acting continuously by isometries on $X$ in a non-elementary way: there is no $G$-fixed point, nor finite orbits on $\bd X$. As explained before, we assume that the action is by type-preserving automorphisms. In this Section, we give an explicit construction of some measurable structures and useful equivariant maps so that the results of Section \ref{section bd theory immeuble} can apply.  
	
	\subsection{Measurable structures on the panel trees}
	Recall that for every $u  \in \bdinf_1$, we denote by $T_u $ the panel tree associated to the vertex at infinity $u$. Endow $T_u$ with the metric $d_u$ given by Equation \ref{eq metric panel tree}. As $X$ is separable, the metric $d_u$ is separable on $T_u$. 
	
	\begin{lem}\label{lem t1 measurable structure}
		The collection of panel trees $T_1 = \bigsqcup_u T_u$ admits a standard Borel structure and the natural projection $p : T_1 \rightarrow \bdinf_1 $ admits a $G$-fiberwise isometric action. 
	\end{lem}
	
	\begin{proof}
		The collection $\{(T_u, d_u)\}_{u \in \bdinf_1}$ is a family of complete, separable metric spaces. Let us construct a fundamental family for $\{(T_u, d_u)\}_{u}$ over $\bdinf_1$. Let $\{y_n\}_n$ be a countable dense set of special vertices of the building $X$. For every $n $, and every $u \in \bdinf_1$, consider 
		\begin{eqnarray}
			x^n_u := \pi_u(y_n) \nonumber
		\end{eqnarray}
		where $\pi_u(y_n)$ is the class (for the strong asymptote relation) of the geodesic ray based at $y$ in the direction of $u$, see Section \ref{section panel trees}. Observe that for every $u $, $\{x^n_u\}_n$ is dense in $T_u$. Fix $n, m \in \mathbb{N}$. Then the function 
		\begin{eqnarray}
			u \in \bdinf_1 \mapsto d_u (x^n_u, x^m_u) \nonumber
		\end{eqnarray}
		is measurable. Indeed, for every $ y  \in X$, the application $\xi \in \bd X \mapsto \gamma_y^\xi$ is continuous for the topology of uniform convergence on compact sets, where $\gamma_y^\xi$ is the unique geodesic ray based at $y$ representing $\xi $. Recall that by definition, 
		\begin{eqnarray}
			d_u (x^n_u, x^m_u) = \underset{s}{\inf} \lim_{t \rightarrow \infty} d(\gamma_{y_n}^u (t+ s), \gamma_{y_m}^u(t)). \nonumber
		\end{eqnarray} 
		Since we consider the infimum of measurable functions, we have that for all $n, m $, the map $u \in \bdinf_1 \mapsto d_u (x^n_u, x^m_u)$ is then measurable. Therefore, the collection $\{x^n\}$ is a fundamental family for the metric field $\{T_u\}_u$ over $\bdinf_1$.
		
		Let us now prove that $G$ acts on the metric field $\{T_u\}_u$. If $u \in \bdinf_1 $ is a vertex at infinity and $\xi \in T_u$ is a class (for the strong asymptote relation) of geodesic rays parallel to $u$, the element $g \in G$ acts as a cocycle $\alpha $ with 
		\begin{eqnarray}
			\alpha (g, u) (u, \xi) = g \cdot(u, \xi) = (gu, g \xi) \nonumber, 
		\end{eqnarray}
		where $g\xi$ is the class (for the strong asymptote relation) of $g \gamma$, for any geodesic ray $\gamma $ in the class of $\xi$. It is clear for all $g \in G$, 
		\begin{eqnarray}
			\alpha (g, u) \in \iso(T_u, T_{gu}) \nonumber. 
		\end{eqnarray}
		It is clear that $\alpha  $ satisfies the cocycle relation: for all $g, g' \in G$, $u \in \bdinf_1$, $\xi \in T_u$, 
		\begin{eqnarray}
			g' \cdot (g \cdot(u, \xi)) = (g'g ) \cdot (u, \xi). \nonumber
		\end{eqnarray} 
		It remains to prove that $G$ acts on this metric field in a measurable way. 
		
		Take $n,  m \in \mathbb{N}$. Observe that for $g \in G$, $u \in \bdinf_1$, 
		\begin{eqnarray}
			g \cdot x^m_{g^{-1}u}&=& g \cdot \pi_{g^{-1}u}(y^m) \nonumber \\
			&=& g \cdot [Q (y^m, g^{-1}u)] \nonumber \\
			&=& [Q(gy^m,u)] \in T_u\nonumber.
		\end{eqnarray}
		As $G$ acts continuously on $X$ and preserves the type of the vertices at infinity, we obtain that 
		\begin{eqnarray}
			(g, u) \mapsto d_u (x^n_u , g \cdot x^m_{g^{-1}u}) \nonumber
		\end{eqnarray}
		is measurable.

		Now by \cite[Lemma 4.11]{duchesne_lecureux_pozetti23}, there is a standard Borel structure on $T_1 := \bigsqcup_{u \in \bdinf_1} T_u$, and a Borel map $p : T_1 \rightarrow \bdinf_1$ such that $p$ admits a $G$-fiberwise isometric action and such that $p^{-1}(u) = T_u$ for every $u \in \bdinf_1$. 
	\end{proof}
	
	Recall that we denote by $\dt_1 $ the union of all the boundaries of the panel trees $T_u, u \in \bdinf_1$. Namely, 
	\begin{equation*}
		\dt_1 = \{ (u, C) \, | \,  u \in \bdinf_1, C \in \partial T_u\}. \nonumber
	\end{equation*}
	
	For $(u ,\xi ) \in T_1$, that is, $u\in \bdinf_1$ and $\xi \in T_u$, we can define the following metric on $\dt_u$:
	\begin{eqnarray}
		d_{u, \xi}(C, D) := \exp(-(C |D)_\xi),\nonumber
	\end{eqnarray} 
	with $(C,D) \in \dt_u$, and where $(. \, | \,. )_\xi$ is the Gromov product on the tree $T_u$ based at $\xi \in T_u$. We define $T_1 \times \dt_1$ as the set 
	\[T_1 \times \dt_1:= \{(u, \xi, C) \mid u \in \bdinf_1, \ \xi \in T_u, \ C \in \partial T_u\}.\]
	
	The next Lemma states that the natural projection $T_1 \times \dt_1 \rightarrow T_1 $ admits a $G$-fiberwise isometric action for this metric.

	\begin{lem}\label{lem dt1 measurable structure}
		The collection $T_1 \times \dt_1 $ admits a standard Borel structure and the natural projection $q : T_1 \times \dt_1 \rightarrow T_1 $ admits a $G$-fiberwise isometric action. 
	\end{lem}
	
	\begin{proof}
		Endow $T_1$ with the measurable structure given by Lemma \ref{lem t1 measurable structure}. Consider the family $\{T_u \times \dt_u\}_{u \in \bdinf_1}$, and call $q$  the projection on the first coordinate 
		\begin{eqnarray}
			q &:& \{T_u \times \dt_u\}_{u \in \bdinf_1} \longrightarrow T_1 \nonumber \\
			&&(u , \xi, C) \longmapsto (u, \xi) \in T_u. \nonumber
		\end{eqnarray}
		For every $u \in \bdinf_1 $ and $\xi \in T_u$, the metric $d_{u, \xi}$ on the fiber $q^{-1}\{(u, \xi)\} = \{(u , \xi)\} \times \dt_u$ is complete and separable. For convenience, we denote the metric spaces $(q^{-1}\{(u, \xi)\}, d_{u, \xi})$ by $\dt_{u, \xi}$. In order to prove that this is a metric field over $T_1$, we need to construct a fundamental family. Consider a countable dense family $\{C_n\}$ of chambers in $\ch(\bdinf)$. For every $u \in \bdinf_1$, define 
		\begin{eqnarray}
			C^n_u := \proj_{u} (C_n) \in \res(u). \nonumber
		\end{eqnarray}
		The projection onto the vertex $u$ is continuous, see \cite[Proposition 3.14]{ciobotaru_muhlerr_rousseau_20}. Moreover, if one fixes $C \in \ch(\bdinf)$, the application $u \in \bdinf_1 \mapsto \proj_{u} (C)$ is continuous. By Proposition \ref{bijection ends chambers}, $C^n_u \in \dt_u$. Consequently, for every $(u, \xi) \in T_1$, define 
		\begin{eqnarray}
			\mathcal{C}^n_{u, \xi} := (u, \xi, C^n_u) \in T_u \times \dt_u \nonumber. 
		\end{eqnarray}
		For all $(u, \xi) $, the collection $\{\mathcal{C}^n_{u, \xi}\}_n$ is dense in the fiber $\dt_{u, \xi}$. Fix $n , m \in \mathbb{N}$, and let $(u , \xi ) \in T_1$. We have that 
		\begin{eqnarray}
			d_{u, \xi}(\mathcal{C}^n_{u, \xi}, \mathcal{C}^m_{u, \xi}) &=& d_{u, \xi}(C^n_u, C^m_u)\nonumber \\
			& = & \exp (-(C^n_u |C^m_u)_\xi) \nonumber \\
			& = & \exp (-(\proj_u(C_n) |\proj_u(C^m))_\xi) \nonumber. 
		\end{eqnarray}
		By composition, the application 
		\begin{eqnarray}
			(u, \xi) \in T_1 \longmapsto d_{u, \xi}(\mathcal{C}^n_{u, \xi}, \mathcal{C}^m_{u, \xi}) \nonumber 
		\end{eqnarray}
		is measurable. The collection $\{\mathcal{C}^n_{u, \xi}\}_n$ is then a fundamental family for the metric field $\{T_u \times \dt_u\}_{u \in \bdinf_1}$ over $T_1 $. 
		
		By Lemma \ref{lem t1 measurable structure}, the projection $T_1 \rightarrow \bdinf_1$ admits a $G$-fiberwise isometric action. For $g \in G$, and $(u, \xi, C) \in \{T_u \times \dt_u\}$, 
		\begin{eqnarray}
			\alpha (g, (u, \xi)) (u, \xi, C) = g \cdot (u, \xi, C) = (gu, g\xi, gC), \nonumber
		\end{eqnarray}
		where $(gu, g\xi) $ was defined in the proof of Lemma \ref{lem t1 measurable structure}, and $gC \in \dt_{gu}$ is the boundary point of $\dt_{gu}$ that corresponds to the chamber $gC \in \res(gu)$ by Proposition \ref{bijection ends chambers}. 
		With this action, for all $g \in G$, 
		\begin{eqnarray}
			\alpha (g, (u, \xi)) \in \iso(\dt_{u, \xi}, \dt_{g \cdot (u, \xi)}). \nonumber
		\end{eqnarray}
		Moreover, the projection $q$ is $G$-equivariant. Last, for every $n, m \in \mathbb{N}$, the measurability of
		\begin{eqnarray}
			(u, \xi) \in T_1 \longmapsto d_{u, \xi}(\mathcal{C}^n_{u, \xi}, g \cdot \mathcal{C}^m_{g^{-1}\cdot(u, \xi)}) \nonumber 
		\end{eqnarray}
		comes from the continuity of the $G$-action on the chambers $\ch(\bdinf)$. This proves that $G$ acts on the metric field $\{T_u \times \dt_u\}_{u \in \bdinf_1}$ with a natural cocycle. By \cite[Lemma 4.11]{duchesne_lecureux_pozetti23}, there is a standard Borel structure on $T_1 \times \dt_1 = \bigsqcup_u T_u \times \dt_u$ such that the projection $T_1 \times \dt_1 \rightarrow T_1 $ is Borel and admits a $G$-fiberwise isometric action. 
	\end{proof}
	
	\subsection{Barycenter map}
	
	Consider the set $(\dt_1)^2 $ defined by 
	\begin{equation*}
		(\dt_1)^2 = \{ (u , C_1, C_2), u \in \bdinf_1, C_1, C_2 \in \partial T_u\}. 
	\end{equation*}
	Denote by $\mathcal{D}(\dt_1)$ the diagonal: 
	$$\mathcal{D}(\dt_1) = \{ (u , C_1, C_2), u \in \bdinf_1, C_1 = C_2 \in \partial T_u\} \subseteq (\dt_1)^2. $$
	Define the application:
	\begin{eqnarray}
		\tau &:& (\dt_1)^2 \to (\dt_1)^2  \nonumber \\
		& & (u , C_1, C_2) \in (\dt_u)^2  \mapsto (u , C_2, C_1) \in (\dt_u)^2 \nonumber. 
	\end{eqnarray}
	Let $(\partial T_1)^{(2)}$ be the measurable set $(\dt_1)^2 -\mathcal{D}(\dt_1)$, and define the space $\partial_2 T_1$ as the set $(\partial T_1)^{(2)}/ \tau $ of unordered pairs of distinct points of $\dt_1$. Similarly, we define by $\partial_n T_1$ as the set of unordered $n$-tuples of points of $\dt_1$ that are pairwise distinct. Endow $\partial_2 T_1$ and $\partial_n T_1$ with the product measurable structure coming from $\dt_1$. 
	
	In the course of the proof, we will need the existence of a well-defined barycenter map on the set $\partial_3 T_1$. 
	
	\begin{prop}\label{prop bary panel trees}
		For all $n\geq3$, there exists a measurable $G$-equivariant map $\partial_n T_1 \to T_1$. 
	\end{prop}
	\begin{proof}
		Let $u \in \bdinf_1$. Let $S\subseteq \partial_n T_u$ be a set of $n$ distinct points in $\partial T_u$. Consider the function 
		$$F_S : x \in T_u \mapsto \sum_{C\neq C' \in S} (C| C')_x .$$
		It is standard to check that the function $F_S$ is proper and convex, see for example \cite[Section 2]{burger_mozes96}. Therefore its minimum $\min(F_S) = \{ x \in T_u \; | \; F_S(x) = \inf F_S\}$ is non-empty, bounded and closed. Since $T_u$ is a tree, we can take the circumcenter of this set, which gives a continuous map
		$$\bary_u : \partial_n T_u \to T_u.$$
		We define
		\begin{eqnarray}
			\bary &:& \partial_n T_1 \to T_1 \nonumber \\
			& & (u, C_1, \dots, C_n) \in\partial_n T_u \mapsto \bary_u(u, C_1, \dots, C_n) \in T_u \nonumber. 
		\end{eqnarray}
		Let $G$ act on $\partial_n T_1$ with the natural diagonal action. Then the map $\bary$ is $G$-equivariant: 
		\begin{eqnarray}
			\bary (g \cdot (u, C_1, \dots, C_n)) &=& \bary_{gu} (g \cdot (u, C_1, \dots, C_n)) \nonumber \\
			& =& g \cdot \bary (u, C_1, \dots, C_n) \nonumber. 
		\end{eqnarray}
		Moreover, the map $\bary$ is measurable for the measurable structures defined by Lemma \ref{lem t1 measurable structure} and \ref{lem dt1 measurable structure}. Indeed, if we take a measurable section $s : \bdinf_1 \to \partial_n T_1 $, then the map 
		$$u \in \bdinf_1 \mapsto \bary_u (s(u))$$
		is the composition of two maps. First, the collection $\{\min(F_{s(u)})\}_{u \in \bdinf_1}$ is a Borel subfield of non-empty closed convex bounded subsets of the metric field $T_1$. Indeed, for any other section $s'$, 
		$$u \in \bdinf_1 \mapsto d_u (s'(u), \min(F_{s(u)}))$$
		is measurable by the definition of $u \mapsto F_{s(u)}$. 
		Next, by Proposition \ref{prop section circ borel}, taking the section of circumcenters gives a Borel map 
		$$\{\min(F_{s(u)})\}_{u \in \bdinf_1} \rightarrow T_u.$$
		This concludes the proof. 
	\end{proof}
	
	Let $Y$ be a proper $\cat(-1) $ space, and let $\prob_3( \bdg Y)$ be the space of positive measures on $\bdg Y$ whose support contains at least three points. In \cite[Proposition 2.1]{burger_mozes96}, Burger and Mozes prove that there exists a $\iso(Y)$-equivariant map 
	$$\bary : \prob_3 (\bdg Y) \to Y$$
	whose restriction to any $\iso (Y)$-orbit is continuous. We shall now prove that there exists a measurable version of this fact, for the metric field $T_1$.

	Let $u \in \bdinf_1$ be a vertex at infinity of type 1 and let $\prob_3( \partial T_u)$ be the space of positive measures on $\partial T_u$ whose support contains at least three points.
	\begin{lem}\label{lem bary proba tu}
		There exists a measurable map $\varphi_u : \prob_3( \partial T_u) \rightarrow \prob(T_u)$. 
	\end{lem}
	
	\begin{proof}
		Let $\nu \in \prob_3( \partial T_u)$. For the moment, assume that $\nu$ is non-atomic. Then consider the probability measure $\tilde{\nu} := \nu \otimes \nu \otimes \nu $ on $(\partial T_u)^3$. Let $\mathcal{D}^3(\partial T_u)$ be the space of triples $(u, \xi_1, \xi_2, \xi_3) \in (\partial T_u)^3$ such that at least two entries coincide. The set $\mathcal{D}^3(\partial T_u) $ is measurable. Since $\nu$ is non-atomic, $\tilde{\nu}$ gives zero measure to $\mathcal{D}^3(\partial T_u) $ and we can see $\tilde{\nu}$ as a measure on $(\partial T_1)^{(3)} $. Passing to the quotient, we can see $\tilde{\nu}$ as a measure on $\partial_3 T_1 $, which we do without changing the notation. Now, thanks to the barycenter map defined in Proposition \ref{prop bary panel trees}, $\varphi_u (\nu) := (\bary_u)_\ast \tilde{\nu}$ is a probability measure on $T_u $.

		Assume now that $\nu$ is purely atomic, that is, for every $\xi \in \supp(\nu)$, $\nu(\xi) > 0$. Since the support of $\nu$ is of cardinal at least $3$, $\nu$ has at least 3 atoms $\xi_1, \xi_2, \xi_3$. Again, let $\tilde{\nu} := \nu \otimes \nu \otimes \nu $ be the product probability measure on $(\partial T_u)^3$. By definition, $\tilde{\nu} (\xi_1, \xi_2, \xi_3) > 0$ and in particular $\tilde{\nu}((\partial T_1)^{(3)} ) >0 $. Recall that taking the restriction of a measure to an open subset is a measurable operation. We can then consider
		$$\bar{\nu} = \frac{1}{\tilde{\nu}((\partial T_1)^{(3)} )} \tilde{\nu}_{|(\partial T_1)^{(3)} },$$
		which is a probability measure on the set of triples of pairwise distinct points in $\dt_u$. Again, we can see $\bar{\nu}$ as a probability measure on $\partial_3 T_1$, and applying the barycenter map $\bary_u$ from Proposition \ref{prop bary panel trees}, one obtains a probability measure $\varphi_u (\nu)$ on $T_u$.

		If now $\nu \in \prob_3( \partial T_u)$ is any probability measure, we can decompose $\nu$ as a sum $\nu = \nu_c + \nu_a$, where $\nu_c$ is non-atomic and $\nu_a$ is purely atomic. As the set $\dt_u$ admits a standard Borel structure, measurability of such a decomposition is classical, see for instance \cite[Theorem~2.12]{dubins_freedman64}. By the previous case, we can assume that the non-atomic part is nonzero $\nu_c > 0$. Therefore, the probability measure $\nu_c / \nu_c (\partial T_u)$ is non-atomic, and by the first case we can measurably associate a probability measure on $T_u$.
	\end{proof}
	We denote by $\prob_3(\partial T_1)$ the subfield of $\prob(\dt_1)$ of measurable sections $u \in \bdinf_1\mapsto \nu_u \in \prob_3(\partial T_u)$ for the measurable structure on $\partial T_1$. 
	\begin{lem}\label{lem bary proba t1}
		There exists a measurable and $G$-equivariant map $\varphi : \prob_3( \partial T_1) \rightarrow \prob(T_1)$ that restricts to $\varphi_u$ on each fiber $\prob_3(\partial T_u)$, for $u \in \bdinf_1$. 
	\end{lem}
	
	\begin{proof}
		On every fiber $u \in \bdinf_1$, the application $\varphi_u : \prob_3( \partial T_u) \rightarrow \prob(T_u)$ is measurable. Moreover, every operation that we did in the proof of Lemma \ref{lem bary proba tu} was $G$-equivariant, i.e. if $\nu = \{\nu_u\}_{u \in \bdinf_1} \in \prob_3(\dt_1)$, then for every $u \in \bdinf_1$, 
		$$g \cdot (\varphi_u (\nu_u)) = \varphi_{gu} (\nu_{gu}).$$
		Measurability of this map comes from Proposition \ref{prop bary panel trees} and the operations in the proof of Lemma \ref{lem bary proba tu}. This proves the measurability and $G$-equivariance of the application. 
	\end{proof}
	
	For self-containment, we now recall the classical construction of a barycenter map for the set of measures on a tree. Let $\varepsilon < \frac{1}{2}$, $\nu \in \prob(T_u)$ and $\xi \in T_u $. For any $R >0$, define 
	$$F^\nu_\xi (R) := \nu(B(\xi, R)) ,$$
	where $B(\xi, R)$ is the ball of radius $R$ and center $\xi$ in $T_u$ for the metric $d_u$ given by equation \eqref{eq metric panel tree}.
	It is clear that $F^\nu_\xi (R) \to 1$ as $R \to \infty $, hence for all $\xi \in T_u $, there exists a well-defined and finite 
	$$R^\nu_{\xi, \varepsilon} := \inf \{ R \; | \; F^\nu_\xi (R) > 1 - \varepsilon\}.$$
	Also, denote by $R^\nu_{\varepsilon} := \inf_{\xi \in T_u} R^\nu_{\xi, \varepsilon}$. We can then define the measurable application 
	$$\tilde{\beta}^u_\varepsilon (\nu):= \overline{\{\xi \in T_u \; | \; R^\nu_{\xi, \varepsilon} < R^\nu_{\varepsilon} +1\}},$$
	It is straightforward to check that $\tilde{\beta}^u_\varepsilon (\nu) $ is a (closed) bounded set in $T_u$. Indeed, if $\xi_1, \xi_2 \in \tilde{\beta}^u_\varepsilon (\nu) $, then by definition 
	$$\nu(\xi \in T_u \; | \; B(\xi_i, R^\nu_{\varepsilon} +1)) > 1/2 \text{ for $i = 1, 2$},$$
	and then $B(\xi_1, R^\nu_{\varepsilon} +1)) \cap B(\xi_2, R^\nu_{\varepsilon} +1)) \neq \emptyset$. Now since trees are $\cat$(0) spaces, Proposition \ref{prop section circ borel} states that we can measurably associate to each bounded set its circumcenter. We then have a measurable application 
	$$\beta^u_\varepsilon : \nu \in \prob(T_u) \mapsto \circum (\tilde{\beta}^u_\varepsilon (\nu) ) \in T_u.$$
	From the definition, $\beta_\varepsilon^u$ is $\iso(T_u)$-equivariant. 
	
	Considering the metric field $\prob(T_1)$ over $\bdinf_1$, what we have obtained is that the map 
	\begin{eqnarray}
		\beta_\varepsilon : \prob(T_1) &\to& T_1 \nonumber \\
		\{\nu_u\}_{u \in \bdinf_1} &\mapsto & \{\beta^u_\varepsilon (\nu_u)\}_{u \in \bdinf_1} \nonumber
	\end{eqnarray}
	is well-defined, measurable, and $G$-equivariant for the natural $G$-action on $\prob(T_1)$. Composing $\beta_\varepsilon $ and $\varphi $ from Lemma \ref{lem bary proba t1}, we have obtained the following result. 
	
	\begin{prop}\label{prop bary proba 3 T1}
		There exists a measurable and $G$-equivariant map 
		\[\prob_3( \partial T_1) \rightarrow T_1\].
	\end{prop}

	We define $T_2$ as the collection of panel trees over the vertices at infinity $\bdinf_2$ of type 2, and we proceed similarly for $\dt_2$ and $\partial_i T_2$. All the results in this section apply for these metric fields in the same way. 
	
\section{Boundary maps in $\tilde{A}_2$-buildings}\label{section bd map immeuble}
	
	The goal of this section is to prove Theorem~\ref{thm intro bd map}. First, we show the existence part:

	\begin{thm}\label{thm bd map immeuble}
		Let $X$ be a  metrically complete and separable building of type $\tilde{A}_2$, $G$ be a locally compact second countable group and $G \curvearrowright X$ a non-elementary action by type-preserving isometries. Let $(B, \nu_B)$ be a $G$-boundary. Then there exists a measurable map $B \rightarrow \ch (\bdinf)$ which is $G$-equivariant. 
	\end{thm}

\subsection{A trichotomy for the boundary map}\label{section bd map trichotomy}
	
	For the rest of this section, we let $X$ be a building of type $\tilde{A}_2$ which is metrically complete and separable, and let $G$ be a locally compact second countable group acting continuously by automorphisms on $X$ in a non-elementary way. Let $(B, \nu_B)$ be a $G$-boundary. 
	
	Recall that the $G$-action is type preserving, so that $G$ acts on $\bdinf_1$ and on $\bdinf_2$. In the proof, we will study the properties of $G$-maps $B \to \bdinf_i$, for $i = 1,2$. The following lemma relies on a result due to B.~Duchesne \cite[Theorem 1.8]{duchesne13}, generalizing Adams-Ballmann Theorem \cite[Theorem 1]{adams_ballmann98} to the context of amenable actions. 
	
	\begin{lem}\label{disjonction}
		Consider the same assumptions as in Theorem \ref{thm bd map immeuble}, and assume that there exists a $G$-map $B \rightarrow \bdinf_1$. Then there is either a $G$-map $B \rightarrow T_1$, or a $G$-map $B \rightarrow \dt_1$, or a $G$-map $B \rightarrow \partial_2 T_1$. 
	\end{lem}
	\begin{proof}
		Let $\phi$ be a $G$-map $\phi: B \rightarrow \bdinf_1$, then by Lemma \ref{lem t1 measurable structure}, the panel trees $\{T_{\phi(b)}\}_{b \in B}$ endowed with their natural metrics form a measurable metric field over $B$. Moreover, there is a natural $G$-action on $T_1$ over $\bdinf_1 $. Each individual $T_u$ is $\cat(0)$, complete and of finite telescopic dimension. Therefore $T_1$ is a Hadamard field of finite telescopic dimension over $\bdinf_1$, and since the $G$-action on $B$ is Zimmer amenable, we can apply \cite[Theorem 1.8]{duchesne13} to the metric field $T_1$ over $\bdinf_1$: either there exists an invariant section of the boundary field $\partial T_1$, or there is an invariant Euclidean subfield of $T_1$. Note that by \cite[Lemma 4.11]{duchesne_lecureux_pozetti23}, an invariant section of the boundary field gives a measurable $G$-map $B \rightarrow \dt_1$.
		
		Otherwise, there is an invariant Euclidean subfield. If this subfield is made of essentially bounded subsets in $T_1$, then taking the section of circumcenters as in Proposition \ref{prop section circ borel} yields a $G$-map $B \rightarrow T_1$. If this subfield consists of lines (1-dimensional Euclidean subspaces), it gives a $G$-map $B \rightarrow \partial_2 T_1$. 
	\end{proof}

	Now our goal is to rule out two of the three previous possibilities. 
	\begin{lem}\label{no T1}
		Under the same assumptions as in Theorem \ref{thm bd map immeuble}, there is no $G$-map $B \rightarrow T_1$. 
	\end{lem}
	
	\begin{proof}
		Let us assume that there is a $G$-map $\psi : B \rightarrow T_1$. Consider the natural projection $p : T_1 \rightarrow \bdinf_1$, and denote $\psi_1 = p \circ \psi : B \rightarrow \bdinf_1$. By Lemma \ref{lem t1 measurable structure}, $p$ is measurable and $G$-equivariant, hence by composition so is $\psi_1$. Take $u,v \in \bdinf_1$ two vertices at infinity of type $1$. Recall that for the boundary of an $\tilde{A}_2 $-building, if $u\neq v$ are both of the same type, the projection $\proj_u (v) $ is a chamber in the residue $\res(u)$. This map is continuous in $v$, and measurable in $u$. Thanks to Proposition \ref{bijection ends chambers}, we can identify the chamber $\proj_u (v) $ to an endpoint of $\dt_u$, which we do without changing the notation. 
		
		Observe that because $\psi_1$ is $G$-equivariant, the set $\{ (b, b') \; | \; \psi_1 (b) = \psi_1(b') \}$ is $G$-invariant. By \cite[Corollary 30]{bader_duchesne_lecureux16}, the action of $G$ on $B \times B$ is ergodic, so this set is either null or conull. If $\nu \otimes \nu$ almost surely, $\psi_1 (b) = \psi_1(b')$, then varying $b$ and $b'$ separately and applying Fubini gives that $\psi_1$ is almost surely constant. But since $\psi_1$ is $G$-equivariant, it means that there is a $G$-fixed point in $\bd X$, which is forbidden by non-elementarity of the $G$-action. As a consequence, $\psi_1(b) \neq \psi_1(b') $ almost surely.
		
		This discussion allows us to define the map:
		\begin{eqnarray}
			\phi: & B\times B &\longrightarrow \dt_1 \nonumber \\
			& (b, b') & \longmapsto (\psi_1(b), \proj_{\psi_1(b)} (\psi_1(b'))) \in \dt_{\psi_1(b)} \nonumber. 
		\end{eqnarray} 
		This application is measurable by composition. As in Lemma \ref{lem t1 measurable structure}, the group $G$ acts on $\dt_1$ via its implicit action on the chambers, using the bijection given by Proposition \ref{bijection ends chambers}. With this action, $\phi$ is $G$-equivariant. Now consider
		\begin{eqnarray}
			\Phi : & B\times B &\longrightarrow T_1 \times \dt_1 \nonumber \\
			& (b, b') & \longmapsto ( \psi(b), \phi(b,b')) \in T_{\psi_1(b)} \times \dt_{\psi_1(b)} \nonumber, 
		\end{eqnarray}
		where  $ T_1 \times\dt_1 $ is the metric field over $T_1$ defined in Lemma \ref{lem dt1 measurable structure}. 
		
		The application $\Phi$ is measurable because $\psi$ and $\phi $ are measurable. The group $G$ acts on $T_1 \times \dt_1$ via the diagonal action $g \cdot (u, \xi, C) = (gu, g\xi, gC)\in T_{gu} \times \dt_{gu}$. Thus, for $g \in G$, and almost every $(b, b') \in B \times B$,
		\begin{eqnarray}
			\Phi(g b , g b')  & = & (\psi (gb ), \proj_{\psi_1(gb)} (\psi_1(gb'))) \nonumber \\
			& = & (g \psi(b ), \proj_{g \psi_1(b)} (g\psi_1 (b'))) \text{ because $\psi $ and $\psi_1$ are $G$-equivariant} \nonumber \\
			& = & (g \psi (b ), g \proj_{\psi_1(b)} (\psi_1(b'))) \nonumber. 
		\end{eqnarray}
		
		As a consequence, $\Phi$ is a $G$-map. We then have the following commutative diagram:
		\[
		\begin{tikzcd}
			B \times B \arrow[r, "\Phi"] \arrow[d, "\pi_B"] 
			& T_1 \times \dt_1 \arrow[d, "\pi_1"] \\
			B \arrow[r, "\psi"]
			& T_1	
		\end{tikzcd}
		\]
		
		By Lemma \ref{lem dt1 measurable structure}, $T_1 \times \dt_1 \rightarrow T_1$ admits a fiberwise isometric $G$-action. But $B$ is a $G$-boundary, hence the projection on the first factor $\pi_B : (b, b') \in B \times B \mapsto b \in B$ is relatively isometrically ergodic. Therefore, there exists an invariant section $s : B \longrightarrow T_1 \times \dt_1$ such that $\nu \otimes \nu $-almost surely, $\Phi (b, b') = s \circ \pi_B(b,b') = s(b)$. 
		\[
		\begin{tikzcd}
			B \times B \arrow[r, "\Phi"] \arrow[d, "\pi_B"] 
			& T_1 \times \dt_1 \arrow[d, "\pi_1"] \\
			B \arrow[r, "\psi"] \arrow[ur, dashed,"s"]
			& T_1	
		\end{tikzcd}
		\]
		
		In particular, the projection $\phi(b,b')= \proj_{\psi_1(b)}(\psi_1(b'))$ does not depend on $b'$. We then have an essentially well defined measurable map $u : B \rightarrow \bdinf_2$ such that for almost every $b \in B$, $u(b) \in \bdinf_2$ is the unique vertex of type 2 belonging to the chamber $\proj_{\psi_1(b)}(\psi_1(b'))$. Then almost surely, $\psi_1(b')$ belongs to a chamber in $\res(u(b))$, see Figure \ref{figure no t1}. 
		
		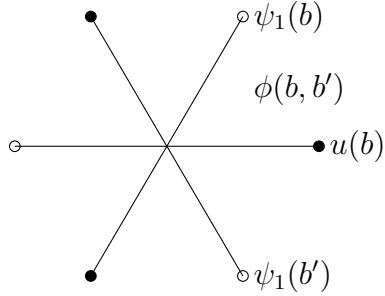
\begin{figure}[h]
			\centering
			\begin{center}
				\begin{tikzpicture}[scale=2]
					\draw (-1,0) -- (1,0)  ;
					\draw (-0.5, -0.86) -- (0.5, 0.86)  ;
					\draw (-0.5, 0.86) -- (0.5, -0.86)  ;
					\draw (0.5, 0.86) circle(1pt) ;
					\draw (0.5, -0.86) circle(1pt) ;
					\draw (-1,0) circle(1pt) ;
					\filldraw (-0.5, 0.86) circle(1pt) ;
					\filldraw (-0.5, -0.86) circle(1pt) ;
					\filldraw (1,0) circle(1pt) ;
					\draw (1,0) node[right]{$u(b)$} ;
					\draw [thick, draw=black, fill=gray!35]
					(1, 0) -- (0,0) -- (0.5, 0.86) ;
					\draw (0.5,0.38) node[right]{$\phi(b,b')$} ;
					\draw (0.5, 0.86) node[right]{$\psi_1(b)$} ;
					\draw (0.5, -0.86) node[right]{$\psi_1(b')$} ;					
				\end{tikzpicture}
			\end{center}
			\caption{An apartment of the building at infinity $\bdinf$, in the proof of Lemma \ref{no T1}.}\label{figure no t1}
		\end{figure}

		In other words, there exists a subset $\Omega_1$ of $B \times B$ of full $\nu \otimes \nu$ measure such that for all $(b,b') \in \Omega_1$, $\psi_1(b')$ belongs to a chamber in $\res(u(b))$. But for all $(b,b') \in \Omega_1$, $\res(u(b))$ is a subset of $\bd X$ of circumradius $\leq \pi/2$. By Fubini, the application $\psi_1 $ then induces an essentially well-defined $G$-map from $B$ to a subset of $\bd X$ of circumradius $\leq \pi/2$. By Proposition \ref{circumcenter}, $G$ fixes a point in $\bd X$, in contradiction with the assumption that the action is non-elementary. 
	\end{proof}

	\begin{lem}\label{no d2T1}
		Under the same assumptions as in Theorem \ref{thm bd map immeuble}, there is no $G$-map $B \rightarrow \partial_2T_1$. 
	\end{lem}
	
	\begin{proof}
		Assume that there exists a measurable $G$-map $\psi : B  \longrightarrow \partial_2 T_1$. Denote by $p : \partial_2 T_1 \rightarrow \bdinf_1$ the natural projection, and $p \circ \psi $ by $\psi_1$. The application $\psi_1 : B \rightarrow \bdinf_1$ is measurable by composition, and $G$-equivariant. 
		
		Again, we define  the application $\phi$ on $B \times B$ by 
		$$\phi : (b,b') \in B\times B \mapsto \proj_{\psi_1(b)} (\psi_1(b')) \in \dt_{\psi_1(b)}.$$
		
		As before, the application $\phi$ is $G$-equivariant, and measurable. Observe that the set $\{(b,b') \in B\times B' \,  | \, \phi(b,b') \in \psi (b)\}$ is $G$-invariant. By double ergodicity of $G \curvearrowright B $ \cite[Corollary 30]{bader_duchesne_lecureux16}, its measure is null or conull. Assume that almost surely, $\phi(b,b') \notin \psi(b)$. Then we get a measurable $G$-map from $B \times B$ to the space of distinct triples of points of $\dt_1$:
		\begin{eqnarray}
			\Psi :&&  	B \times B  \longrightarrow  \partial_3 T_1  \nonumber \\
			&&(b,b')  \longmapsto  (\psi(b), \phi(b,b')) \in \partial_3 T_{\psi_1(b)}. \nonumber
		\end{eqnarray}
		Now consider the barycenter map	
		\begin{eqnarray}
			\bary : && \partial_3 T_1  \longrightarrow T_1 \nonumber \\
			&&(C_1, C_2, C_3) \in \partial_3 T_u  \longmapsto  \bary_u(C_1, C_2, C_3) \in T_u \nonumber, 
		\end{eqnarray}
		By Proposition \ref{prop bary panel trees}, the barycenter map just defined is measurable and $G$-equivariant, hence $\bary \circ \Psi $ is a measurable $G$-map. It implies that the following diagram
		\[
		\begin{tikzcd}
			B \times B \arrow[r, "\bary \circ \Psi"] \arrow[d, "\pi_B"] 
			& T_1 \arrow[d, "\pi_1"] \\
			B \arrow[r, "\psi_1"]
			& \bdinf_1	
		\end{tikzcd}
		\]
		is commutative. But $T_1 \rightarrow \bdinf_1$ admits a fiberwise isometric $G$-action by Lemma \ref{lem t1 measurable structure}. Therefore, by relative isometric ergodicity of the Poisson boundary, there exists an invariant section $s : B \longrightarrow T_1$, which contradicts Lemma \ref{no T1}.
		
		As a consequence, $\phi(b,b') \in\psi(b)$ almost surely, so that $\phi(b,b')$ essentially does not depend on $b'$. But the image of $\phi(b,b')$ lies in $\res(\psi_1(b))$, which is a subset of $\bd X $ of circumradius $\leq \pi/2$. Now by the same arguments as for Lemma \ref{no T1}, we get a $G$-fixed point in the boundary and a contradiction. 
	\end{proof}
	
	It is clear that Lemma \ref{no T1} and Lemma \ref{no d2T1} remain true for the vertices of type 2: there is no $G$-map $B \rightarrow T_2$ nor $G$-maps $B \rightarrow \partial_2 T_2$. 
	We can now finish the proof of Theorem \ref{thm bd map immeuble}. 
	
	\begin{proof}[Proof of Theorem \ref{thm bd map immeuble}]
		By assumption, $G$ acts continuously without invariant flats on the affine building $X$, and $X$ is a complete $\cat$(0) space of finite telescopic dimension. We can then apply Theorem \ref{furst map}: there exists a measurable $G$-map $\psi : B \rightarrow \bd X$. 
		
		The partition of the visual boundary $\bd X $ between regular points (classes of geodesic rays whose endpoints belong to the interiors of the chambers at infinity), and singular points (classes of geodesic rays that represent vertices at infinity) is measurable because $\bdinf_1$ and $\bdinf_2$ are closed in $\bd X$. In other words, we have the measurable decomposition $\bd X =  \bdinf_1 \sqcup \bdinf_2 \sqcup (\bd X - \bdinf_1 \cup \bdinf_2 )$, and the set of regular points $\bd X - \bdinf_1 \cup \bdinf_2 $ projects equivariantly and measurably onto $\ch(\bdinf)$. 
		
		Observe that the set $\{b \in B \; | \; \psi (b) \in \bd X - (\bdinf_1 \sqcup \bdinf_2)\}$ is measurable, and $G$-invariant. By ergodicity of $G \curvearrowright B $, its measure is null or conull. If for almost every $b \in B$, $\phi (b) $ is a regular point, then we get a $G$-map $\psi : B \rightarrow \ch (\bdinf)$, which is what we want. 
		
		Otherwise, by ergodicity of the action $G \curvearrowright B$, we get either a map $ \psi : B \rightarrow \bdinf_1$ or $ \psi : B \rightarrow \bdinf_2$. Assume the former, the argument being identical in both cases. 
		
		We can apply Lemma \ref{disjonction}: there is either a $G$-map $B \rightarrow T_1$, or a $G$-map $B \rightarrow \dt_1$, or a $G$-map $B \rightarrow \partial_2 T_1$. Along with Lemmas \ref{no T1} and \ref{no d2T1}, there is a $G$-equivariant measurable map $B \rightarrow \dt_1$. By Proposition \ref{bijection ends chambers}, the chambers at infinity in the residue of a vertex $u \in \bdinf_1$ are $G$-equivariantly and measurably associated with $\dt_u$, hence $\psi$ gives a $G$-equivariant measurable map $ B \rightarrow \ch (\bdinf)$. 
	\end{proof}

\subsection{Uniqueness of the boundary map}
	Keep the same assumptions and notations as in Theorem \ref{thm bd map immeuble}. 
	
	\begin{thm}\label{thm uniqueness chamber map}
		The $G$-map $\psi : B \longrightarrow \ch(\bdinf)$ given by Theorem \ref{thm bd map immeuble} is essentially the unique measurable $G$-equivariant boundary map $ B \longrightarrow \ch(\bdinf)$. 
	\end{thm}
	\begin{proof}
		Recall that for each chamber $C \in \ch({\bdinf})$, $C$ contains a unique vertex of type 1 and a unique vertex of type 2. Hence the boundary map $\psi : B \longrightarrow \ch(\bdinf)$ induces two $G$-maps $\psi_1 : B \longrightarrow \bdinf_1$ and $\psi_2 : B \longrightarrow \bdinf_2$. It is then sufficient to prove that those are the only $G$-maps $ B \longrightarrow \bdinf_1$ and $ B \longrightarrow \bdinf_2 $. The argument being the same in both cases, we only prove that $\psi_1: B \longrightarrow \bdinf_1$ is unique. Let us consider $\psi'_1 : B \to \bdinf_1$ be a $G$-equivariant measurable map. By ergodicity of the $G$-action on $B$, the set 
		$$\{ b \in B \mid \psi_1(b) = \psi'_1(b)\}$$ 
		is of measure 0 or 1. Let us assume by contradiction that this set is of measure 0.  For (almost every) $b \in B$, let us denote by $u(b) \in \bdinf_2$ the unique vertex of type $2$ that is adjacent to both $\psi(b) $ and $\psi'(b)$ in any apartment containing them. Alternatively, $u(b)$ can be defined as the unique vertex of type 2 adjacent to the chamber $\proj_{\psi_1(b)}(\psi'_1(b))$ as in the proof of Lemma \ref{no T1}. As we assume that almost surely $\psi_1(b) \neq \psi'_1(b)$ we obtain that almost surely, 
		$$ \proj_{u(b)}(\psi_1(b))\neq \proj_{u(b)}(\psi'_1(b)).$$
		Using the identification between the residue $\res(u(b)) $ and the ends of the panel tree $T_{u(b)}$, we get a measurable $G$-equivariant map 
		\begin{eqnarray}
			\Phi &:& B \longrightarrow \partial_2T_2 \nonumber \\
			&& b \mapsto (u(b), \proj_{u(b)}(\psi_1(b)), \proj_{u(b)}(\psi'_1(b))) \in \partial_2T_{u(b)}. \nonumber
		\end{eqnarray}
		In view of Lemma \ref{no d2T1}, we get a contradiction. 
	\end{proof}

	In the course of the proof, we have shown the following result of independent interest.
	\begin{cor}
		With the same assumptions as in Theorem \ref{thm bd map immeuble}, then the map $\psi_1 = \proj_i \circ \psi$ is essentially the unique $G$-equivariant measurable map $\psi_i : B \longrightarrow X^\infty_i$, for $i= 1,2$. 
	\end{cor}
	With the notations of the Section \ref{section intro alg}, it means that any $\Gamma$-map $B \to G/Q$, for $Q$ a parabolic subgroup contained in the minimal parabolic subgroup $P$, comes from the composition the map $B \to G/P $ given by Theorem \ref{thm bd map immeuble} and the projection $G/ P \to G/Q$ coming from the inclusion $P < Q$. 
	
	\section{Stationary measure on the chambers at infinity}\label{section uniq stat meas A2}
	\subsection{Uniqueness of the stationary measure}
	
	Let $X$ be a separable complete building of type $\tilde{A}_2$ and let $G$ be a locally compact second countable group acting continuously by isometries on $X$ in a non-elementary way. 
	The main result of this section is the following. 
	
	\begin{thm}\label{thm uniqueness stationary measure immeuble}
		Let $\mu $ be an admissible symmetric measure on $G$. Then there is a unique $\mu$-stationary measure $\nu$ on  $\ch(\bdinf)$. Moreover, we have the decomposition 
		$$\nu = \int_{b \in B} \delta_{\psi(b)} d\nu_B(b),$$
		where $(B, \nu_B)$is the Poisson-Furstenberg boundary of $\mu$ and $\psi$ is the unique boundary map given by Theorem \ref{thm uniqueness chamber map}. In particular, $(\chinf, \nu)$ is a $(G, \mu)$-boundary. 
	\end{thm}
	
	Before proving this theorem, recall the following classical result, which goes back to H.~Furstenberg \cite[Lemma 1.33]{furstenberg73}.

	\begin{lem}\label{furst}
		Let $B(G, \mu)$ be the Poisson-Furstenberg boundary associated to the probability measure $\mu$ on $G$, and let $M$ be a Borel $G$-space. Then for any $\mu$-stationary measure $\nu \in \prob(M)$, there is a $G$-equivariant map $ \phi : B \rightarrow \prob(M)$ such that $\nu = \int \phi (b) d \nu(b)$. Conversely, if $\phi : B \rightarrow \prob(M)$ is measurable and $G$-equivariant, then the measure $$\nu = \int \phi (b) d \mathbb{P}(b)$$ is $\mu$-stationary. 
	\end{lem}

	\begin{proof}[Proof of Theorem \ref{thm uniqueness stationary measure immeuble}]
		By Lemma \ref{furst}, it is enough to prove that there is a unique $G$-equivariant map $ B \longrightarrow \prob(\ch(\bdinf))$. Let $\psi : B \longrightarrow \ch (\bdinf)$ be the measurable $G$-equivariant map given by Theorem \ref{thm bd map immeuble}. Observe that the map 
		$$b \in B \mapsto \delta_{\psi(b)} \in \prob(\ch(\bdinf))$$
		is measurable and $G$-equivariant. Let $\phi : B \rightarrow \prob(\ch(\bdinf))$ be any $G$-equivariant map. The goal is to show that almost surely, $\phi $ is given by $\phi(b) = \delta_{\psi(b)}$. 
		
		Recall that each chamber in $\bdinf$ contains exactly one vertex of type 1 and one vertex of type 2, so denote by $\pi_1 : \ch(\bdinf) \longrightarrow \bdinf_1 $ and $\pi_2 : \ch(\bdinf) \rightarrow \bdinf_2 $ the corresponding maps. Let $\psi_i: B \rightarrow \bdinf_i$ be defined by $\psi_i= \pi_i \circ \psi : B \rightarrow \bdinf_i$, for $i = 1, 2 $. By composition, $\psi_i$ is a measurable, $G$-equivariant map. Consider also the map 
		$$\phi_1 : b \in B \mapsto  ({\proj_{\psi_1(b)}})_\ast \phi(b) \in \prob(\ch(\bdinf)).$$
		The map $\phi_1$ is a measurable $G$-equivariant map between $B$ and the probability measures on the chambers, but now the support of $\phi_1(b)$ is contained in $\res({\psi_1(b)})$. Now identify $\res(\psi_1(b))$ with the ends $\dt_{\psi_1(b)}$ of the panel tree at $\psi_1(b)$ given by Proposition \ref{bijection ends chambers}. Therefore, $\phi_1$ can be written as 
		\begin{eqnarray}
			\phi_1 : B &\rightarrow &\prob(\dt_1) \nonumber \\
			b & \mapsto &\phi_1(b) \in \prob(\dt_{\psi_1(b)}). \nonumber 
		\end{eqnarray} 
		By ergodicity of the Poisson boundary, the cardinal $k(b)$ of the support of $\phi_1(b)$ is almost surely constant. If almost surely $k(b) \geq 3$, then the barycenter map constructed in Proposition \ref{prop bary proba 3 T1} associates to the support of $\phi(b)$ a canonical point in $T_{\psi_1(b)}$. This gives a measurable $G$-equivariant map $B \longrightarrow T_1$, which is impossible due to Lemma \ref{no T1}.

		If the support of the measure $\phi_1(b)$ is almost surely of cardinal $2$, then we have a $G$-map $B \longrightarrow \partial_2 T_1$, which is in contradiction with Lemma \ref{no d2T1}.

		Then $\phi_1(b)$ has to be a Dirac mass, which we denote by $\delta_{\xi(\omega) }$ for $(\psi_1 (\omega), \xi (\omega )) \in \dt_{\psi_1(\omega)}$. Again, by Proposition \ref{bijection ends chambers}, this gives a $G$-map $B \longrightarrow \ch(\bdinf)$. But by Theorem \ref{thm uniqueness chamber map}, any such map is unique and equal to $\psi$. Therefore, $\phi_1(b)$ is the Dirac mass at $\psi(b)$. 
		
		What we obtained so far is that for any $C \in \ch(\bdinf)$, 
		\begin{eqnarray}\label{eq uniq stat measure A2}
			\phi_1(b) (\proj_{\psi_1(b)} (C)) &=& \delta_{\psi(b)} (\proj_{\psi_1(b)} (C)).
		\end{eqnarray}
		By definition, $ \phi_1(b)  = ({\proj_{\psi_1(b)}})_\ast \phi(b) $. Then by Equation \eqref{eq uniq stat measure A2}, the support of $\phi(b)$ is contained in the set of chambers $C \in \ch(\bdinf)$ such that $ \proj_{\psi_1(b)} (C) = \psi(b) $. 
		
		By repeating the same argument as before for $\psi_2 : B \rightarrow \bdinf_2$, we obtain that $\phi(b) $ must also be supported on the set of chambers $C \in \ch(\bdinf)$ such that $ \proj_{\psi_2(b)} (C) = \psi(b) $. Combining these two results, we get that $\phi(b)$-almost surely, $\proj_{\psi_1(b)} (C) = \proj_{\psi_2(b)} (C) = \psi(b)$, hence $\phi(b)$ is supported on the chamber $\psi(b)$, meaning that $\phi(b)$ is essentially equal to $\delta_{\psi(b)}$. This proves the uniqueness of $\psi$, therefore of the stationary measure. 
	\end{proof}

	Recall that if we denote by $(B, \nu_B)$ the Poisson-Furstenberg boundary associated to $(G,\mu)$, we have a canonical measurable map 
	$$\bnd : \Omega \rightarrow B$$
	such that $\bnd_\ast \mathbb{P} = \nu_B$. 
	An immediate application of Theorem \ref{thm lim meas furst} gives the following useful result.  
	\begin{cor}\label{cor furstenberg lim meas}
		Keep the same notations as in Theorem \ref{thm uniqueness stationary measure immeuble}. Let $\nu $ be the unique $\mu$-stationary measure on $\chinf$ associated to $\mu$, and let $\psi : B \to \ch(\bdinf)$ be the unique boundary map given by Theorem \ref{thm uniqueness chamber map}. Then $\mathbb{P}$-almost surely, $Z_n (\omega)\nu \rightarrow \delta_{\psi(\bnd(\omega))} $ in the weak-$\ast$ topology. 
	\end{cor}

	\subsection{Opposite chambers at infinity}
	Recall that two chambers $C, C'\in \chinf$ are opposite if they are at maximal gallery distance, equivalently if they belong to a unique apartment of $\bdinf$. Since apartments are convex, the convex hull of two chambers is contained in any apartment containing both of them. 
	
	\begin{prop}\label{prop bd map antipodal}
		Let $B$ be a a $G$-boundary. Then the measurable equivariant map $\psi: B \to \chinf$ given by Theorem \ref{thm uniqueness chamber map} is almost surely antipodal: for almost every $b,b' \in B$, $\phi(b)$ and $\phi(b')$ are opposite. 
	\end{prop}
	\begin{proof}
		Consider the unique $G$-equivariant map $\psi : B \rightarrow \ch(\bdinf)$ given by Theorem \ref{thm bd map immeuble}, and the following measurable map : 
		\begin{eqnarray}
			\phi : & B\times B &\longrightarrow W_0 \nonumber \\
			& (b, b') & \longmapsto \delta(\psi(b), \psi(b')) \nonumber, 
		\end{eqnarray}
		where $\delta$ is the Weyl distance function associated to the spherical Coxeter group $(W_0, S)$ of type $A_2 $, associated to the spherical building at infinity $\bdinf$. Then by $G$-equivariance of $\psi$, the measurable map $\phi$ is $G$-invariant. By double ergodicity of the Poisson boundary \cite[Corollary 30]{bader_duchesne_lecureux16}, $\phi $ is then essentially constant. 
		
		If its essential value is $e$, then almost surely, $\psi(b) = \psi(b')$, and there is a $G$-fixed chamber on the boundary, which is impossible because the action is non-elementary. 
		
		Let us denote by $l_S$ the word metric on $W$ given by the set $S$. From now on, let $B_0, B_1 \subseteq B$ be conull sets such that for all $(b,b') \in B_0 \times B_1$, the essential value of $\phi$ is $\phi(b,b') $. If almost surely, $l_S(\phi(b,b'))= 1 $, then almost surely $\psi(b)$ and $\psi(b')$ are adjacent. Since $\phi$ is constant, then almost surely $\psi(b)$ and $\psi(b')$ contain a vertex of constant type, say of type $1$. Let $b \in B_0$, denote by $u(b) \in \psi(b)$ the unique vertex of type $1$. Then for all $b' \in B_1$, $\psi(b') \in \res(u(b))$. Again, $\res(u(b))$ has radius $\leq \pi/2$ so by Proposition \ref{circumcenter}, we obtain a contradiction because the action is non-elementary. 
		
		If $l_S(\phi)$ is almost surely equal to 2. Then the gallery between $\psi(b)$ and $\psi(b')$ is of constant type $(1,2)$ or $(2,1)$. But this is impossible because $\delta(\psi(b'), \psi(b))= \delta(\psi(b), \psi(b'))^{-1}$, so by ergodicity, we have that $\delta(\psi(b), \psi(b'))= \delta(\psi(b), \psi(b'))^{-1}$. 
		
		Consequently, $\psi(b)$ and $\psi(b')$ are almost surely opposite. 
	\end{proof}
	
	Denote by $\nu$ the unique $\mu$-stationary measure on $\ch(\bdinf)$ given by Theorem \ref{thm uniqueness stationary measure immeuble}. The following is an important feature of the stationary measure $\nu$. 
	\begin{cor}\label{cor proba opposite chambers}
		Let $\nu$ be the unique $\mu$-stationary measure on $\ch(\bdinf)$ associated to $\mu$, as given by Theorem \ref{thm uniqueness stationary measure immeuble}. Then $\nu \otimes \nu$-almost every pair of chambers in $\ch(\bdinf)$ are opposite. 
	\end{cor}
	
	\begin{proof}
		Let $B$ be the Poisson-Furstenberg boundary associated to $(G, \mu)$. The proof is a direct combination of Proposition \ref{prop bd map antipodal} applied to $B$ and the decomposition 
		$$\nu = \int_{b \in B} \delta_{\psi(b)} d\nu_B(b)$$
		given by Theorem \ref{thm uniqueness stationary measure immeuble}. 
	\end{proof}
	
	\section{Convergence to a chamber at infinity}\label{section conv}
	\subsection{Regular sequences in buildings}\label{section reg}
	
	In this section, $X$ is any affine building. We present the notion of regular sequences and relate it with Lyapunov regularity. Regular sequences first appeared in \cite{kaimanovich89} for symmetric spaces, and were then developed in \cite{parkinson_woess15} for affine buildings. 
	
	We fix $\mathfrak{a}^+$ a fundamental closed Weyl chamber for the affine reflection group $(W,V)$. For any two points $x,y \in X$, the type $\theta(x,y)$ is a vector in $\mathfrak{a}^+$. 
	\begin{Def}
		Let $\lambda \in \mathfrak{a}^+$ be a fixed vector. An embedding $\gamma : \R_+ \to X$ is called a \emph{$\lambda$-regular ray} if for any $s,t \in \R_+$, 
		$$\theta(\gamma(t), \gamma(t+s)) = s\lambda \in \mathfrak{a}^+. $$
		A sequence $(x_n) $ in $X$ is called \emph{$\lambda$-regular} if there exists a $\lambda$-regular ray $\gamma$ such that 
		$$ d(x_n, \gamma(n)) = o(n ).$$
	\end{Def}
	
	The following theorem is the main result of \cite{parkinson_woess15}, where the authors only consider discrete buildings. In the context of symmetric spaces, it was proven by Kaimanovich in \cite[Theorems 2.1 and 2.4]{kaimanovich89}
	
	\begin{thm}\label{thm regular Lyapunov}
		Let $(x_n)$ be a sequence in an affine building $X$, and let $\lambda \in \mathfrak{a}^+$. Then the following statements are equivalent: 		
		\begin{enumerate}
			\item \label{reg implies lyap} the sequence $(x_n) $ is $\lambda$-regular;
			\item \label{lyap implies reg} $d(x_n, x_{n+1})= o(n) $ and 
			$$\frac{1}{n} \theta (o, x_n) \longrightarrow \lambda.$$
		\end{enumerate}  
	\end{thm}
	
	Most of the proof in \cite[Theorem 3.2]{parkinson_woess15} goes through to non-discrete buildings. We denote by $\{\alpha_i\}_i$ the set of roots, and by $\{\lambda_i\}_{i \in I_0}$ be the basis of $V$ dual of $\{\alpha_i\}_i$, so that $\langle \alpha_i , \lambda_j \rangle = \delta_{i,j}$, which we call \emph{coweights}. The ``only if'' part is straightforward.
	\begin{proof}
		(\ref{reg implies lyap} $\Rightarrow $ \ref{lyap implies reg}) Let us first assume that $(x_n) $ is $\lambda$-regular. Let $\gamma : [0, \infty[ \rightarrow X $ be a $\lambda$-ray. If $\lambda = 0$, the result is immediate. We can then assume without loss of generality that $\lambda \neq 0$ and $o = \gamma (0)$. Then by the triangular inequality, 
		$$d(x_n, x_{n+1}) \leq d(x_n, \gamma(n))+d(\gamma(n), \gamma(n+1))+d(\gamma(n+1), x_{n+1}) = o(n).$$
		Since $\gamma$ is a $\lambda$-ray starting at $o$, we have that $\theta (o, \gamma(n)) = n \lambda$ for all $n$. 
		Now $$\theta (o, x_n) = \sum_i \langle \theta (o, x_n), \alpha_i\rangle \lambda_i.$$
		By the law of cosines in $\cat(0)$ spaces (see for instance \cite[Exercise II.1.9]{bridson_haefliger99}), 
		$$ d(x_n, \gamma(n))^2 \geq  d(o, x_n)^2 + d(o, \gamma(n))^2 - 2 d(o, x_n)d(o, \gamma(n)) \cos(\angle_o (x_n, \gamma(n))).$$
		Since $d(x_n, \gamma(n))= o(n)$ and $d(o, \gamma(n)) = n \|\lambda\|$, we obtain that $\angle_o (x_n, \gamma(n)) \rightarrow_n 0 $. As a consequence, for all $i \in I_0$ we have that
		$$ \langle \theta (o, x_n), \alpha_i\rangle = \langle \theta (o, \gamma(n)), \alpha_i\rangle + o (n), $$
		whence the result.

	(\ref{lyap implies reg} $\Rightarrow $ \ref{reg implies lyap}) The approach of this part given in \cite{parkinson_woess15} goes through without modifications, as it mostly relies on $\cat(0)$ geometry and Lemma \ref{lem park woess} below, valid for any affine building, and whose proof in the non-discrete case is exactly the same (see also \cite[2.4.11]{rousseau23}).  
	\end{proof}
	
	\begin{lem}\label{lem park woess}
		Let $o, x,y \in X$ be such that $\theta(o, x) = \theta (o, y) = \lambda$ and for which $[o, x] \cap [o, y] = \{o\}$. Then there is $C>0 $ depending only on the direction of $\lambda $ (not on its length) such that $d(x, y) \geq C d(o, x)$. Actually,
		$$C = \min \{ 2\sin (\angle(\lambda, w \lambda)) \mid w \in W_0 - \{1\}\}.$$
	\end{lem}
	
	\subsection{Convergence of the random walk}
	
	Let now $X$ be a separable complete building of type $\tilde{A}_2$ and let $G$ be a discrete countable group acting continuously by isometries on $X$ in a non-elementary way. Let $\mu$ be an admissible measure on $G$. As in Section \ref{section bd theory}, we let $(\Omega, \mathbb{P}) $ be the probability space $(G^{\mathbb{N}}, \mu^{\otimes \mathbb{N}})$ with the product $\sigma$-algebra, and $(Z_n(\omega))$ be the random walk on $G$ associated to $\mu$. 
	
	We recall that $\mu$ is said to have finite $n$-th moment for the action $G \curvearrowright (X,d)$ if for some (equivalently, any) $o \in X$, 
	$$\int_G d(o, go)^n d\mu(g) < \infty.$$
	
	The \emph{escape rate} (or \emph{drift}) of the random walk is defined as 
	$$l_X(\mu) :=\lim_{n\rightarrow \infty} \frac{1}{n}\int_\Omega d( o,Z_n(\omega)o) d \mathbb{P}(\omega)$$
	if $\mu$ has finite first moment, and $l_X(\mu)= \infty $ otherwise. 
	
	We are going to use the following theorem from H. Izeki. 
	\begin{thm}[{\cite[Theorem A]{izeki23}}]\label{thm izeki}
		Let $Y$ be a complete CAT(0) space which is either proper or of finite
		telescopic dimension, and let $G$ be a discrete countable group equipped with a symmetric and admissible probability measure $\mu$ with finite second moment. Assume that $G \curvearrowright Y $ is a non-elementary action. Then the drift of the random walk is positive:
		$$l_Y(\mu) > 0 .$$
	\end{thm}
	
	In our case, we consider a metrically complete $\adt$-building, hence of finite telescopic dimension (equal to 2), and the theorem applies. We are ready to complete the proof of Theorem~\ref{thm intro rw}, which we restate for convenience. 
	
	\begin{thm}[Simplicity of the Lyapunov spectrum]\label{thm simplicity}
		Let $X$ be a separable complete building of type $\tilde{A}_2$ and let $G$ be a discrete countable group acting continuously by isometries and non-elementarily on $X$. Let $\mu$ be an admissible measure on $G$ with finite second moment. Then there exists a regular $\lambda \in \mathfrak{a}^{++}$ such that for almost every $\omega \in \Omega$, $(Z_n(\omega))$ is $\lambda$-regular. In particular, $(Z_n o)$ converges almost surely to a regular point of the visual boundary. 
	\end{thm}

	\begin{proof}
		By Theorem \ref{thm izeki}, the drift $l_\mu := l_X(\mu)$ of the random walk $(Z_n(\omega)o)$ is positive. We can then apply \cite[Theorem 2.1]{karlsson_margulis99}: for almost every $\omega$, there exists a geodesic ray $\gamma^\omega$ such that 
		$$\lim \frac{1}{n} d(\gamma^\omega(l_\mu n), Z_n(\omega)o) = 0.$$
		In other words, if we denote $\lambda^\omega= \theta(\gamma^\omega(0),\gamma^\omega(l_\mu)) \in \mathfrak{a}^+$ whenever it is defined, then $(Z_n(\omega)o)$ is almost surely $\lambda^\omega$-regular. 
		In particular, the random walk converges almost surely to a point of the visual boundary. We then obtain a $G$-equivariant measurable map 
		$$ F: B \to \bd X, $$
		which to almost every $b \in B$ associates the limit point of $(Z_n(\omega))$ for $\omega \in \bnd^{-1}(b)$. By ergodicity of the $G$-action on the Poisson boundary, this boundary point is either almost surely a regular point, a vertex at infinity of type 1 or a vertex at infinity of type 2. Note that this decomposition corresponds exactly to the type of $\lambda_\mu$: $\lambda_\mu$ is regular (resp. of type 1,2) if and only if $(Z_n(\omega))$ converges to a regular point of the boundary (resp. to a vertex of type 1,2).

		Now consider $\check{\mu}= i_\ast \mu$ for $i(g) = g^{-1}$. The measure $\mu$ being symmetric by assumption, we have $\lambda_\mu = \lambda_{\check{\mu}}$. However, if we denote by $\iota : \mathfrak{a}^+ \to \mathfrak{a}^+$ the involution, we have that $\lambda_{\check{\mu}} = \iota(\lambda_\mu)$, see for instance \cite[Corollary 9.11]{benoist_quint}. In $\tilde{A}_2$ buildings, the involution changes the type. In particular, the Lyapunov vector $\lambda_\mu$, being involution-invariant, is regular. 
	\end{proof}
	\begin{rem}
		We emphasize that the geometric property of $\tilde{A}_2$-buildings that is crucial for the simplicity of the Lyapunov spectrum is only used at the very end of the proof: there exists no proper facet of the closed fundamental Weyl chamber that is left invariant by the opposition involution. In the $\tilde{C}_2$ and $\tilde{G}_2$ cases, this is false, and showing the simplicity of the Lyapunov spectrum requires more assumptions. 
	\end{rem}
	
	\subsection{Combinatorial convergence}
	We now prove that, under the same assumptions, the random walk converges to a chamber at infinity in the polyhedral sense too. For a detailed account on this notion, we refer to \cite{caprace_lecureux11} (for discrete buildings) and to \cite[\S 3.4]{rousseau23}. In what follows, we denote by $\Sigma_o X$ the residue building based at $o$, as defined in \cite[Corollary 1.11]{parreau00}. It is a spherical building, and for every special vertex $o \in X$, we have a canonical morphism of simplicial complexes 
	$$\Sigma_o : \bdinf \to \Sigma_o X$$
	sending any facet at infinity $F^\infty $ to $\germ_o(Q(o,F^\infty))$. 
	
	\begin{cor}
		We keep the assumption of Theorem \ref{thm simplicity}. Let $o,x \in X$. Then for almost every $\omega$, there exists $n_0$ such that for all $n \geq n_0$, the projection $\Sigma_o(Z_n(\omega)x)$ of $Z_n(\omega)x$ on the residue building $\Sigma_o X$ is the chamber $\Sigma_o (\tilde{\psi}(\omega))$, where $\tilde{\psi}(\omega) = \psi(\bnd(\omega)) \in \chinf$ is the boundary map given by Theorem \ref{thm bd map immeuble}. 
	\end{cor}
	
	\begin{proof}
		By Theorem \ref{thm simplicity}, there exists $\lambda \in \mathfrak{a}^{++} $ a regular vector and $\Omega' \subseteq \Omega $ of full measure satisfying: for all $\omega \in \Omega'$, there exists a $\lambda$-ray $\gamma^\omega$ such that the random walk $(Z_n(\omega)x)$ sublinearly tracks $\gamma^\omega$. Fix  $\omega \in \Omega'$, and write $\gamma := \gamma^\omega$. We can assume without loss of generality that $\gamma(0) = o$, so that $\gamma$ is contained in the interior of the Weyl chamber $S= Q(o, C_\omega)$, where $C_\omega = \tilde{\psi}(\omega)$. The vector $\lambda$ is regular, hence the distance between $\gamma(n)$ and any of the walls delimiting $S$ grows linearly. Since $d(x_n, \gamma(n)) = o(n)$, there exists $N$ such that for all $n \geq N$, $d(x_n, \gamma(n))$ is strictly less than the distance from $\gamma(n)$ and any of these walls. Let $n \geq N$. We need to prove that the germ of the geodesic segment $[o, x_n]$ is contained in $\germ (S)$. Let $S_n$ be a Weyl chamber based at $o$ containing the $[o, x_n]$. By \cite[Proposition 1.15]{parreau00}, there exists an apartment $A_n$ containing $S$ and the germ of $S_n$ in $o$. Now consider the retraction $\rho_{A_n, S}$ onto $A_n$ centered on $S$. The retraction does not increase the distance 
		$$ d(\rho_{A_n, S}(x_n), \rho_{A_n, S}(\gamma(n)))= d(\rho_{A_n, S}(x_n), \gamma(n))\leq d(x_n, \gamma(n)),$$
		hence for $n \geq N$, the point $\rho_{A_n, S}(x_n) $ is contained in the interior of $S$. As $\rho_{A_n, S}$ induces a retraction of $\Sigma_o X$ onto $\Sigma_o A_n$, we obtain that the germ of $[o, x_n]$ is strictly contained in the  $\germ_o(S)$ and the result follows.
	\end{proof}
	
	\subsection{Stationary measures on $X \cup \ch(\bdinf)$}
	In this subsection, we show that the measure on $\ch(\bdinf)$ given by Theorem \ref{thm uniqueness stationary measure immeuble} is actually the only stationary measure on the bordification $X\cup \ch(\bdinf)$.

	Let us then endow $X \cup \chinf$ with the metrizable separable Hausdorff topology on $X \cup \chinf$ described in Section \ref{section intro immeuble}.
	
	\begin{prop}
		We keep the notations of Theorem \ref{thm simplicity}. Let $\nu$ be any $\mu$-stationary probability measure supported on the bordification $X \cup \chinf$ for the cone topology. Then the support of $\nu $ is in $\chinf$.
	\end{prop}
	\begin{proof}
		Let $o \in X$. Up to decomposing the measure $\nu$, we assume by contradiction that $\nu$ is supported in $X$. Then there exists $r > 0$ such that the open ball $O :=B(o,r)$ satisfies $\nu (O) \geq \frac{2}{3}$. Since $\nu$ is $\mu$-stationary, we can apply Theorem \ref{thm lim meas furst}: there exists an essentially well-defined $G$-equivariant measurable map $\omega \mapsto \nu_\omega$ such that $Z_n (\omega)\nu \rightarrow \nu_\omega $ almost surely in the weak-$\ast$ topology, and such that we have the decomposition $$\nu = \int_{\Omega} \nu_\omega d \mathbb{P}(\omega). $$
		As by Theorem \ref{thm simplicity}, $(Z_n o)$ converges almost surely to a regular point in a chamber at infinity, for almost every $\omega$, there exists $n_0 (\omega)$ such that for all $n \geq n_0$, we have $Z_n (\omega) O \cap O = \emptyset$ (take for instance $n_0 $ such that  $d(o, Z_no ) \geq 3r$ for all $n \geq n_0$). In particular, for almost every $\omega \in \Omega $ and for all $n \geq n_o (\omega)$, we have $\nu (Z_n^{-1}(\omega)O)=Z_n (\omega)\nu (O) < \frac{1}{3}$. As a consequence, $\nu_\omega (O) \leq \frac{1}{3}$ for almost every $\omega$. Since $\nu$ decomposes as $\int_{\Omega} \nu_\omega d \mathbb{P}(\omega)$, we get that $\nu(O)\leq \frac{1}{3}$, a contradiction. 
	\end{proof}
	In view of Theorem \ref{thm uniqueness stationary measure immeuble}, we derive the following. 
	\begin{cor}\label{cor unique stat measure bordif}
		We keep the notations of Theorem \ref{thm simplicity}. There is a unique $\mu$-stationary probability measure $\nu$ on $X \cup \chinf$, and $\nu $ is given by Theorem \ref{thm uniqueness stationary measure immeuble}. 
	\end{cor}

	\bibliographystyle{alpha}
	\bibliography{bibliography}
\end{document}